\RequirePackage{fix-cm}
\documentclass[smallextended, numbook]{svjour3}       
\smartqed  
\usepackage{graphicx}
\usepackage{algorithm}
\usepackage{amsmath}
\usepackage{amsfonts}
\usepackage{graphicx}
\usepackage{epstopdf}
\usepackage{algorithmic}
\usepackage[super]{nth}
\usepackage{color,colortbl}
\usepackage{soul}
\usepackage{pgfplotstable}
\usepackage{cleveref}
\usepackage{booktabs}
\usepackage{multicol}
\usepackage{ifplatform}
\usepackage[toc,acronym,nonumberlist]{glossaries}

\usepackage{xargs}                      
\usepackage[colorinlistoftodos,prependcaption,textsize=tiny]{todonotes}
\newcommandx{\unsure}[2][1=]{\todo[linecolor=red,backgroundcolor=red!25,bordercolor=red,#1]{#2}}
\newcommandx{\change}[2][1=]{\todo[linecolor=blue,backgroundcolor=blue!25,bordercolor=blue,#1]{#2}}
\newcommandx{\info}[2][1=]{\todo[linecolor=OliveGreen,backgroundcolor=OliveGreen!25,bordercolor=OliveGreen,#1]{#2}}
\newcommandx{\improvement}[2][1=]{\todo[linecolor=Plum,backgroundcolor=Plum!25,bordercolor=Plum,#1]{#2}}
\newcommandx{\thiswillnotshow}[2][1=]{\todo[disable,#1]{#2}}

\newcommand{\convexhull}[1]{Conv(#1)}
\newcommand{\bracket}{\mathcal{X}}
\newcommand{\bracketstep}[1]{\sigma_{#1}}
\newcommand{\bracketupdate}[1]{U_{#1}}
\newcommand{\bracketinnerdiam}[1]{b(#1)}
\newcommand{\bracketdiam}[1]{b(#1)}
\newcommand{\bracketouterdiam}[1]{diam(#1)}

\newcommand{\altbracket}{\pmb t}

\newcommand{\bracketcheck}{\varphi}

\newcommand{\swappoints}{X}
\newcommand{\swapupdates}{W}
\newcommand{\bracketset}{\mathcal{B}}
\newcommand{\semibracketset}{\tilde{\mathcal{B}}}
\newcommand{\upmconstant}{\alpha}
\newcommand{\upmscaling}{h}

\newcommand{\supmtag}{S}
\newcommand{\eupmtag}{E}
\newcommand{\dupmtag}{D}
\newcommand{\eupmsequence}{I}
\newcommand{\eupmsequencesetgen}{\mathcal{I}_n}
\newcommand{\eupmsequencesetpart}{\mathcal{I}_5}
\newcommand{\eupmsequenceset}[1]{\mathcal{I}_{#1}}
\newcommand{\eupmsubstring}{J}
\newcommand{\eupmsubstringset}{\mathcal{S}}
\newcommand{\argmin}{\operatornamewithlimits{argmin}}
\newcommand{\condensed}{1}

\newcommand{\extbracketmiddle}{x^M}
\newcommand{\extbracketleft}[2]{x^L_{#1,#2}}
\newcommand{\extbracketleftone}{x^L_1}
\newcommand{\extbracketlefttwo}{x^L_2}
\newcommand{\extbracketleftthree}{x^L_3}
\newcommand{\extbracketlefti}{x^L_\idxi}
\newcommand{\extbracketleftj}{x^L_\idxj}
\newcommand{\extbracketright}[2]{x^L_{#1,#2}}
\newcommand{\extbracketrightone}{x^R_1}
\newcommand{\extbracketrighttwo}{x^R_2}
\newcommand{\extbracketrightthree}{x^R_3}
\newcommand{\extbracketrighti}{x^R_\idxi}
\newcommand{\extbracketrightj}{x^R_\idxj}

\newcommand{\altbracketleftone}{q}
\newcommand{\altbracketlefttwo}{p}
\newcommand{\altbracketrightone}{r}
\newcommand{\altbracketrighttwo}{s}

\newcommand{\idxi}{i}
\newcommand{\idxj}{j}

\newacronym{upm}{UPM}{Underestimating Polynomial Method}
\newacronym{dupm}{DUPM}{Dynamic Underestimating Polynomial Method}
\newacronym{supm}{SUPM}{Static Underestimating Polynomial Method}
\newacronym{eupm}{EUPM}{Extremal Underestimating Polynomial Method}
\newacronym{sqp}{SQP}{Sequential Quadratic Programming}
\newacronym{sqpgs}{SQP-GS}{Sequential Quadratic Programming with Gradient Sampling}

\crefname{problem}{problem}{problems}%
\crefname{Problem}{Problem}{Problems}%

\newcommand{\nice}{locally unimodal }

\newcommand{\mindist}{\delta}
\newcommand{\truesol}{x^*}
\newcommand{\funcprojection}{Q}

\newcommand{\algmiff}{the Mifflin-Strodiot method }
\begin{document}

\title{A seven-point algorithm for piecewise smooth univariate minimization\thanks{This publication is based on work supported by the EPSRC Centre for Doctoral
		Training in Industrially Focused Mathematical Modelling (EP/L015803/1) in collaboration with Siemens.}
}


\author{Jonathan Grant-Peters         \and
        Raphael Hauser 
}


\institute{J. Grant-Peters \at
              Mathematical Institute, University of Oxford, United Kingdom \\
              \email{peters@maths.ox.ac.uk}           
           \and
           R. Hauser \at
              Mathematical Institute, University of Oxford, United Kingdom
}

\date{Received: date / Accepted: date}

\maketitle

\begin{abstract}
In this paper, we construct an algorithm for minimising piecewise smooth functions for which derivative information is not available. The algorithm constructs a pair of quadratic functions, one on each side of the point with smallest known function value, and selects the intersection of these quadratics as the next test point. This algorithm relies on the quadratic function underestimating the true function within a specific range, which is accomplished using a adjustment term that is modified as the algorithm progresses.
\end{abstract}

\section{Introduction}

\label{sec:intro}

	\ifnum \condensed=0
		In this chapter, we derive an algorithm suitable for  minimising the piecewise smooth univariate function $f$ in the interval $[a,b]$. 
		As we will later use this algorithm as a line search for the multivariate problem, we need it to be first robust and then fast. 
		We begin by exploring existing algorithms and assess whether a suitable one exists.
	\else
		When solving the problem of minimizing a function of many variables, many existing solutions belong to a class of algorithms called line-search algorithms \cite{nocedal2006numerical,fletcher2013practical,gould2006introduction}
		which reduce the problem by iteratively restricting the search to 1-dimensional subspaces determined by a search direction to determine a step size.
		\ifnum \condensed=1
			This 1-dimensional problem may be solved either exactly or inexactly.
			Most methods rely on the latter approach, as it is computationally cheaper while no less effective.
			
			While exact line search is rarely used, due to its cost, there exist certain special cases where it is more effective such as Linear Programming. 
			Moreover, it has been observed by Yu et al \cite{yu2010quasi} that in the context of non-smooth optimisation, using exact line-search may result in stepping to a location from where a better subdifferential approximation may be constructed, thus resulting in selecting better search directions at future iterations.
		\else
			There are two ways to approach this one dimensional sub-problem: the exact and inexact approaches. 
			For inexact line-search, we merely look in the search direction for a point at which the objective function is smaller than at the current iterate and satisfies certain other conditions \cite{wolfe1969convergence,wolfe1971convergence}.
			The strength of this approach is that fewer function evaluations are required for the one-dimensional problem, while the inexactness of the line-search in practice makes little difference.
			
			The other approach is the exact line-search, which involves rigorously solving a one dimensional sub-problem with each iteration. 
			While this approach is usually more computationally expensive than inexact line-search, there exists certain cases where it is superior, such as for linear programming problems.
			Furthermore, it has been observed by Yu et al \cite{yu2010quasi} that in the context of non-smooth optimisation, using exact line-search may result in selecting better search directions at future iterations.
		\fi

		In this paper, motivated by designing an exact line search algorithm for non-smooth optimisation, we examine the effectiveness of existing univariate optimisation algorithms on a class of non-smooth objective functions and introduce a new algorithm which is specifically designed for this class.
		The particular problem for which we develop this method is a black box, piecewise smooth objective function for which derivatives are not available.
	\fi
	
	When conducting an exact line-search, there are two distinct steps. 
	First of all we must  construct a \textit{bracket}, which means finding a closed interval which is guaranteed to contain a local minimum of the univariate objective function.
	The second is find the local minimum within this interval.
	In this paper, we focus mainly on the second step which we rigorously state below.
	
	\begin{problem}
		\label{prob:1dopt}	
		Given is a function $f:\mathbb{R}\rightarrow \mathbb{R}$ of the form $f = \max_{i = 1, \ldots, k} f_i(x)$ where $f_i:\mathbb{R} \rightarrow \mathbb{R}$ are smooth functions, an oracle for calling function values of $f(x)$, and an interval $[a,b]$. Design an optimisation method for finding a local minimiser of $f$ in $[a,b]$ which makes use of only oracle calls for $f$.
	\end{problem}

	One class of algorithms which solves \Cref{prob:1dopt} is that of 1d global optimisers.
	This class includes using interval analysis based methods such as the Moore-Skelboe algorithm \cite{moore1966interval}.
	However, these methods require that the objective function is known, and are not applicable to black box functions.
	
	Another category of solutions is that of global Lipschitzian methods such as those summarised by Hansen et al \cite{hansen1992global}. 
	While these might be feasible if we could guess a suitable Lipschitz Constant, in practice there is little gained from them.
	Even if more than one local optimum existed, we have no reason to expect that finding one minimum as opposed to another would make an exact line search based method more effective.
	Therefore, while we will make use of some of the ideas behind such algorithms, we reject them as they involve needless extra work.

	
	When choosing a local univariate solver, one faces a trade off between speed and robustness. 
	The basic methods are Golden Section \cite{kiefer1953sequential,gill2019practical} which converges Q-linearly for any continuous function, and  interpolating methods \cite{gill2019practical,jarratt1967iterative} which, when successful, converge super linearly.
	The latter's stability  depends on their ability to construct a polynomial which approximates the objective function well locally. 
	If the approximation does not fit the objective function well, then these algorithms may fail. 
	
	In practice, the most effective methods are bracketing interpolation hybrid methods such as Brent's Method \cite{brent1976new,brentalgorithms} and Hager's cubic method \cite{hager1989derivative}  which combine the speed of interpolation methods with the stability of Golden Section by making use of a fall back option when the interpolation is not working. 
	While convergence is guaranteed for such hybrid methods, they will still need their interpolations to match the objective function well if super linear convergence is to be obtained.
	Therefore, we can't expect them to converge quickly for an objective function of the form stated in \Cref{prob:1dopt}.
	
	There exist other bracketing-interpolation hybrid methods for non differentiable functions of the form $f(x) = \max_i f_i(x)$ \cite{murray1979steplength,yu2010quasi}.
	However, these algorithms assume that we can compute each $f_i$ (from which $f$ is computed) separately. 
	Therefore, these are unsuitable for \Cref{prob:1dopt} given that the oracle is defined to only return the value of $f(x)$.
	
	The only algorithm which seems to be  optimised for our setting is the five point method of Mifflin and Strodiot \cite{mifflin1993rapidly}, to which we will from now on refer to as ``the Mifflin-Strodiot method''.
	This algorithm is designed to converge rapidly to the local minimum $x^*$ even when the objective function is non-differentiable at $x^*$. 
	
	The remainder of this paper is structured as follows: In \Cref{sec:bracketing}, we define the notion of a \textit{bracket} rigorously and describe in brief the existing methods which are currently most applicable to our problem.
	\Cref{sec:upm,sec:eupm,sec:dupm} focus on the derivation of a new univariate optimiser, for which we present three variations.
	We compare our new methods against relevant competitors in \Cref{sec:univariate_numerics}.


\section{Bracketing Methods}
\label{sec:bracketing}
We begin by rigorously defining a bracket:
\begin{definition}[Bracket]
	\label{def:bracket}
	Let $f:\mathbb{R}\rightarrow \mathbb{R}$ be a function, and $x^L, x^M, x^R \in \mathbb{R}$. We call the trio of points $x^L, x^M, x^R$ a \textbf{bracket} of $f$ if they satisfy the following:
	\begin{multicols}{2}
		\begin{enumerate}
			\item $x^L < x^M < x^R$,
			\item $f(x^L) \geq f(x^M) \leq f(x^R)$.
		\end{enumerate}
	\end{multicols}
\end{definition}
The significance of a bracket is that \Cref{def:bracket} guarantees that there exists a local minimum of $f$ in $(x^L, x^R)$. 
\begin{definition}[Set of Brackets]
	\label{def:bracketset}
	We define $\bracketset_f$ to be the set of all brackets for the function $f$.
\end{definition}

In this paper, we refer to the elements of a bracket $\bracket$ as $x^L$, $x^M$ and $x^R$. 
If a bracket $\bracket_i$ is associated with a particular iteration of an algorithm, then we write its elements as $x^L_i$, $x^M_i$ and $x^R_i$. 
Given this notation, we define the function $\bracketinnerdiam{\bracket} = x^R - x^L$.

\begin{algorithm}
	\caption{General Bracketing Method}
	\label{alg:bracketing_method_template}
	\begin{algorithmic}
		\REQUIRE $f$ Objective function, $\bracketstep{}$ Step Function, $\bracketupdate{}$ Update Function,\\ $\bracket_0 $  Initial Bracket, $\epsilon>0$ Tolerance.
		
		\STATE $i = 0$;
		
		\WHILE{$\bracketinnerdiam{\bracket_i}> 2\epsilon$\label{alg:line:bracket_termination}} 
		
		\STATE $\tilde x_i = \bracketstep{}(\bracket_i)$.
		
		\STATE $\bracket_{i+1} = U(\tilde x_i,\bracket_i)$. \label{alg:line:bracket_update}
		
		\STATE $i = i+1;$
		\ENDWHILE
		
		\RETURN $\bracket_i$.
	\end{algorithmic}
\end{algorithm}

Bracketing methods are the set of algorithms which are vaguely in the form of \Cref{alg:bracketing_method_template}, given the inputs of a step function $\bracketstep{}$ and an update function $\bracketupdate{}$.
The step function should select a new point within $(x^L, x^M)\cup(x^M, x^R)$, while the update function should return a bracket.
When discussing existing bracketing methods and constructing our new one,  we use the structure of \Cref{alg:bracketing_method_template} where one algorithm is distinguished from another based on how the step function $\bracketstep{}$ and update function $\bracketupdate{}$ are defined.

The purpose of this paper is to introduce a new bracketing method  called the Underestimating Polynomial Method (UPM) for  \Cref{prob:1dopt} for which, to our knowledge,  there does not currently   exist a robust and fast solution.
When assessing the effectiveness of this algorithm, we will compare it to a small selection of existing bracketing methods including Golden Section, Brent's Method \cite{brent1976new,brentalgorithms}, and the Mifflin-Strodiot method \cite{mifflin1990superlinear}.

Of these algorithms, Golden Section may be considered to most robust as it is guaranteed to converge Q-linearly with a rate of approximately $0.618$. 
Meanwhile Brent's method is the most effective for smooth functions converging super-linearly, while sometimes performing surprisingly well for non-smooth functions.
Finally, Mifflin's method is theoretically the most comparable to ours, in that it is equipped for non-smooth functions.
However, we will see that it lacks the robustness of the previous two algorithms.

	The idea of bracketing methods such as Brent's method \cite{brentalgorithms,brent1976new},  and Hager's Cubic method \cite{hager1989derivative} is that the polynomial which interpolates the bracket points approximates the objective function well locally.
Therefore the local minimum of the quadratic is a sensible location to evaluate next.  
The error of this approximation can be quantified and bounded, yielding super-linear convergence guarantees \cite[Theorem 4.1]{brentalgorithms}.
While this analysis works well for smooth functions, our problem is a piecewise smooth, black box function.
\begin{definition}
	\label{def:pwsmooth}
	A function $f:  \mathbb{R}^n\supseteq \Omega \rightarrow \mathbb{R}$ is called \emph{piecewise smooth} if $f$ is continuous on $\Omega$ and there exists a disjoint finite family of sets $\{\Omega_i;\;i=1,2,\ldots,k\}$ such that $f$ is smooth on $\Omega_i$ for all $i \in \mathcal{I}$ and $\cup_{i=1}^k \Omega_i$ is dense in $\Omega$.
	We call any point $x$ such that $x\in \overline{\Omega_i} \cap \overline{\Omega{j}}$ for some $i,j\in \mathcal{I}$ a \emph{kink}.
\end{definition}
The theory behind algorithms such as Brent's method collapses due to the existence of kinks, because interpolating across a kink has no meaning and yields no meaningful error bound. 
Therefore, if the local minimum of a function $f$ is a kink, we expect algorithms like Brent's method to converge to it slowly.

\ifnum \condensed=0
	\subsection{Underestimating Polynomial Method}
\else
	\section{Static Underestimating Polynomial Method}
\fi
\label{sec:upm}

\ifnum \condensed=0
	The idea of bracketing methods such as Brent's method \cite{brentalgorithms,brent1976new},  and Hager's Cubic method \cite{hager1989derivative} is that the polynomial which interpolates the bracket points approximates the objective function well locally.
	Therefore the local minimum of the quadratic is a sensible location to evaluate next.  
	The error of this approximation can be quantified and bounded, yielding super-linear convergence guarantees \cite[Theorem 4.1]{brentalgorithms}.
	While this analysis works well for smooth functions, our problem is a piecewise smooth, black box function.
	\begin{definition}
		\label{def:pwsmooth}
		A function $f:  \mathbb{R}^n\supseteq \Omega \rightarrow \mathbb{R}$ is called \emph{piecewise smooth} if $f$ is continuous on $\Omega$ and there exists a disjoint finite family of sets $\{\Omega_i;\;i=1,2,\ldots,k\}$ such that $f$ is smooth on $\Omega_i$ for all $i \in \mathcal{I}$ and $\cup_{i=1}^k \Omega_i$ is dense in $\Omega$.
		We call any point $x$ such that $x\in \overline{\Omega_i} \cap \overline{\Omega{j}}$ for some $i,j\in \mathcal{I}$ a \emph{kink}.
	\end{definition}
	The theory behind algorithms such as Brent's method collapses due to the existence of kinks, because interpolating across a kink has no meaning and yields no meaningful error bound. 
	Therefore, if the local minimum of a function $f$ is a kink, we expect algorithms like Brent's method to converge to it slowly.
\fi

The premise of the \acrshort{upm} is inspired by the approach used in \algmiff \cite{mifflin1990superlinear}.
We approximate the objective function with two polynomials, one on each side of the bracket.
If the local minimum is a kink, then the combination of these polynomials may be valid approximations and useful for locating the minimum. 
To construct these polynomials, we need more than the three points contained in a bracket and therefore extend the definition of a bracket to include seven points:
\begin{definition}[Extended Bracket]
	\label{def:extended_bracket}
	Given the function $f$, we call $\bracket \in \mathbb{R}^7$ an extended bracket written in the following form
	\begin{equation}
	\label{eq:bracket_extension}
	\bracket = \left( \begin{array}{ccccccc}\extbracketleftthree& \extbracketlefttwo&\extbracketleftone&\extbracketmiddle&\extbracketrightone&\extbracketrighttwo&\extbracketrightthree \end{array} \right)^T,
	\end{equation} 
	if the following conditions apply:
	\begin{subequations}
		\begin{align}
		\extbracketleftthree< \extbracketlefttwo<\extbracketleftone<&\extbracketmiddle<\extbracketrightone<\extbracketrighttwo<\extbracketrightthree\label{eq:x_condition_1}\\
		f(\extbracketleftone)\geq f(&\extbracketmiddle)\leq f(\extbracketrightone).\label{eq:x_condition_2}
		\end{align}
	\end{subequations}
\end{definition}
\begin{remark}
	If $\bracket$ is an extended bracket, then  $(\extbracketleftone, \extbracketmiddle, \extbracketrightone)$ form a bracket. 
\end{remark}
For expressing an extended bracket's size, we define the following functions.
\begin{definition}
	[Extended Bracket Length]\label{def:extended_bracket_length}
	Let $\bracket $ be an extended bracket of $f$. 
	We define: $\bracketouterdiam{\bracket} = diam(\convexhull{\bracket})$, and $\bracketinnerdiam{\bracket} = \extbracketrightone - \extbracketleftone$. 
\end{definition}
\begin{remark}
	\Cref{def:extended_bracket_length} is the natural extension of the function $\bracketinnerdiam{\bracket}$ from \Cref{sec:bracketing} where $\bracket$ is merely a bracket as opposed to an extended bracket.
\end{remark}

%
%

\ifnum \condensed=0  

The subsections below are structured in the following way. In \Cref{sec:upm_core}, we give a vague outline of the method, in addition to introducing the notation we use in later sections. Afterwards in \Cref{sec:supm,sec:eupm,sec:dupm} we present three variations of the \acrshort{upm}, the order of which is motivated by the derivation process.
\else
\fi

\ifnum \condensed=0
	\subsection{\acrshort*{upm} core method}
	\label{sec:upm_core}
\fi

%

The \acrshort{upm} mostly conforms to the form of  \Cref{alg:bracketing_method_template}. 
The main difference, apart from replacing the bracket with an extended bracket, is the fact that the \acrshort{upm} also depends on an input parameter $\upmconstant$. 
This parameter is used by the step function $\bracketstep{}$ when constructing the next point to evaluate. 
The three variations of the \acrshort{upm} presented in this paper differ in how they tread $\upmconstant$. 

The \acrfull*{supm} requires an initial value of $\upmconstant$ to be supplied by the user, which is then kept constant for the duration of the algorithm.
For the remainder of this section, we denote the step function $\bracketstep{}$ and update function $\bracketupdate{}$ which apply to the \acrshort{supm} by $\bracketstep{\supmtag}$ and $\bracketupdate{\supmtag}$ respectively.

Given the objective function $f$, our strategy is to use the points $\extbracketleftone, \extbracketlefttwo, \extbracketleftthree$ and $\extbracketrightone, \extbracketrighttwo,\extbracketrightthree$  to construct the model functions $q^L(x;\bracket, \upmconstant)$, and $q^R(x;\bracket, \upmconstant)$ which approximate $f_L$ and $f_R$ respectively. 
We will employ Newtons Divided Difference notation which we summarize below.

\begin{definition}
	\label{def:first_divided_difference}
	Given $a,b, c \in \mathbb{R}$, and $f:\mathbb{R} \rightarrow \mathbb{R}$, the \nth{1} and \nth{2} divided differences are defined by:
	\begin{align*}
	f[a,b] &= \frac{f(a) - f(b)}{a - b},&&\text{and}&& f[a,b,c] = \frac{f[a,b]- f[a,c]}{b - c}.
	\end{align*}
\end{definition}

The reason we use this notation is that the \nth{1} and \nth{2} divided differences are natural approximations for the \nth{1} and \nth{2} derivatives of $f$ respectively.
We will state the result in a later section when more rigour is needed.

Now we define our model functions $q^L$ and $q^R$ as the following:
\begin{equation}
\label{eq:UPM_model_function_definitions}
q^k(x; \bracket,\upmconstant) = f(x^k_{1}) + f[x^k_{1}, x^k_{2}](x-x^k_{1}) + (f[x^k_{1}, x^k_{2}, x^k_{3}]-\upmconstant \upmscaling(\bracket))(x-x^k_{1})(x - x^k_{2}),
\end{equation}
where $k \in \{L,R\}$, $\upmconstant>0$ is a constant and $\upmscaling(\bracket)$ is a scaling function with the property that $\upmscaling(\bracket) \rightarrow 0$ as $\bracketouterdiam{\bracket} \rightarrow 0$. 
Note that $q^k$ has the form of the \nth{2} order Newton Interpolating Polynomial combined with the adjustment term $\upmconstant \upmscaling(\bracket)$.
Finally we express $\bracketstep{\supmtag}$ in terms of $q^L$ and $q^R$:

\begin{equation}
\label{eq:UPM_gp}
 \bracketstep{\supmtag}(\bracket;\upmconstant) = \argmin \limits_{x \in [\extbracketleftone, \extbracketrightone]} \max(q^L(x;\bracket, \upmconstant), q^R(x;\bracket, \upmconstant) ).
\end{equation}


\begin{remark}
	The model functions $q^L$ and $q^R$ are designed to underestimate $f_L$ and $f_R$ within $[\extbracketleftone, \extbracketrightone]$. 
	We will show how we achieve this in \Cref{lem:alpha_threshold}.
\end{remark}
\begin{figure}[h!]
	\includegraphics[width=1\textwidth]{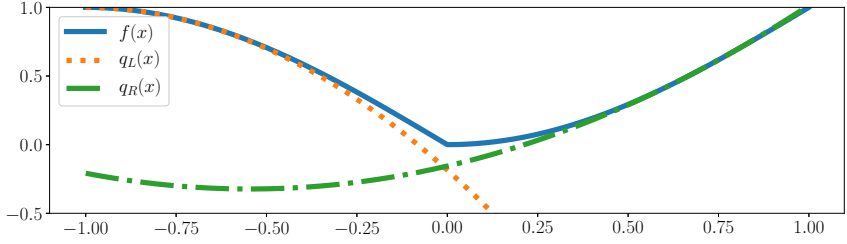}
	\caption{We plot the objective function in addition to the quadratic function approximating the left and right. These two quadratics intersect at $x=0.005\ldots$.}
	\label{fig:idea_of_method}
\end{figure}

\begin{example}
	Consider the function $f(x) = \max \left(\sin (\tfrac{1}{2}x \pi), 1 - \cos (\tfrac{1}{2} x \pi)\right)$.
	This function is  piecewise smooth  and unimodal in the interval $[-1,1]$ with a local minimum at $0$. 
	At this local minimum, $f$ is not differentiable.
	
	Now suppose that $\{\extbracketlefti\}_{i=1}^{3} = \{-0.75, -0.9, -1\}$ and $\{\extbracketrighti\}_{i=1}^{3} = \{0.6, $ $0.8,$ $ 0.95\}$. 
	We construct $q^L$ and $q^R$ by interpolating $\{(x^k_i, f(x^k_i))\}_{i=1}^{3}$ for $k=L,R$. 
	These two quadratics intersect at $x\approx 0.005$ as shown in \Cref{fig:idea_of_method}.
	This value is taken to be the next point at which we evaluate $f$. 
\end{example}

All that remains is for us to define the function $\bracketupdate{\supmtag}$. 
\ifnum \condensed = 0 
With each iteration there are exactly 4 different ways in which we might update $\bracket_n$. The one we choose depends on the values of $\tilde x = \bracketstep{\supmtag}(\bracket;\upmconstant)$ and $f(\tilde x)$ compared to $f(\extbracketmiddle)$. In particular we define:
\fi
\begin{equation}
\label{eq:UPM_iteration_update}
 \bracketupdate{\supmtag}(\tilde x; \bracket, \upmconstant) := \left\{ \begin{array}{lcc}   
 \bracketupdate{1}(\tilde x;\bracket) & \text{if}&\tilde x<\extbracketmiddle\;\text{and}\;f(\tilde x) < f(\extbracketmiddle)\\
 \bracketupdate{2}(\tilde x;\bracket) & \text{if}& \tilde x>\extbracketmiddle\;\text{and}\;f(\tilde x) < f(\extbracketmiddle)\\
 \bracketupdate{3}(\tilde x;\bracket) & \text{if}& \tilde x>\extbracketmiddle\;\text{and}\;f(\tilde x) > f(\extbracketmiddle)\\	
 \bracketupdate{4}(\tilde x;\bracket) & \text{if}& \tilde x<\extbracketmiddle\;\text{and}\;f(\tilde x) > f(\extbracketmiddle)
\end{array}    \right.,
\end{equation}  
where
\begin{subequations}
	\label{eq:UPM_update_functions}
	\begin{align}
	\bracketupdate{1}(\tilde x;\bracket) :=\left( \begin{array}{ccccccc}\extbracketleftthree& \extbracketlefttwo&\extbracketleftone&\tilde x&\extbracketmiddle&\extbracketrightone&\extbracketrighttwo \end{array} \right)^T ,\\	
	\bracketupdate{2}(\tilde x;\bracket) :=\left( \begin{array}{ccccccc} \extbracketlefttwo&\extbracketleftone&\extbracketmiddle&\tilde x&\extbracketrightone&\extbracketrighttwo&\extbracketrightthree \end{array} \right)^T ,\\
	\bracketupdate{3}(\tilde x;\bracket) :=\left( \begin{array}{ccccccc}\extbracketleftthree& \extbracketlefttwo&\extbracketleftone&\extbracketmiddle&\tilde x&\extbracketrightone&\extbracketrighttwo \end{array} \right)^T ,\\
	\bracketupdate{4}(\tilde x;\bracket) := \left( \begin{array}{ccccccc} \extbracketlefttwo&\extbracketleftone&\tilde x&\extbracketmiddle&\extbracketrightone&\extbracketrighttwo&\extbracketrightthree \end{array} \right)^T.
	\end{align}
\end{subequations}
\begin{remark}
	The update function $\bracketstep{\supmtag}(\bracket)$ is the natural extension of the one used by \algmiff \cite{mifflin1993rapidly} to 7 points. 
\end{remark}

\begin{lemma}
	\label{lem:upm_maintains_bracket}
	Let  $\bracket_1$ be an extended bracket. 
	If $\tilde x = \bracketstep{\supmtag}(\bracket_1) \in (\extbracketleftone,\extbracketmiddle)\cup(\extbracketmiddle, \extbracketrightone)$, then  $\bracket_2 = \bracketupdate{\supmtag}(\tilde x; \bracket_1, \upmconstant)$ is an extended bracket and $\bracketinnerdiam{\bracket_2}< \bracketinnerdiam{\bracket_1}$.   
\end{lemma}
\ifnum \condensed=0
\begin{proof}
	
	If $\tilde x \in (\extbracketleftone , \extbracketmiddle)$, then the update function $\bracketupdate{}$ applies either $\bracketupdate{1}$ or $\bracketupdate{4}$ to $\bracket_1$.
	In both cases, $\tilde x$ is placed between $\extbracketleftone$ and $\extbracketmiddle$ which implies that $\bracket_2$ satisfies \Cref{eq:x_condition_1}.

	Given that $\bracket_1$ satisfies \Cref{eq:x_condition_2}, it follows that $\extbracketmiddle$ is the location within $[\extbracketleftone, \extbracketrightone]$ where the smallest value of $f$ has been evaluated. 
	If $f(\tilde x) < f(\extbracketmiddle)$, then $\bracketupdate{\supmtag}$ applies $\bracketupdate{1}$, which implies $\bracket_2$ satisfies \Cref{eq:x_condition_2}. 
	Alternatively if $f(\tilde x) > f(\extbracketmiddle)$, then $\bracketupdate{\supmtag}$ applies $\bracketupdate{4}$, in which case $\bracket_2$ still satisfies \Cref{eq:x_condition_2}.
	Since $\bracket_2$ satisfies \Cref{eq:x_condition_1,eq:x_condition_2}, it is an extended bracket.
	
	The equivalent argument may be applied to the case where $\tilde x \in (\extbracketmiddle , \extbracketrightone)$, exchanging $\bracketupdate{1}$ and $\bracketupdate{4}$ for $\bracketupdate{2}$ and $\bracketupdate{3}$ respectively.  
\end{proof}
\fi

%

\ifnum \condensed=0
	\subsection{\acrlong*{supm}}
	\label{sec:supm}
	
	The first variation of the \acrshort{upm} which we discuss  is the \acrfull*{supm}. 
	Here the word static refers specifically to the treatment of the parameter $\upmconstant$. 
	The \acrshort{supm} requires a user defined value of $\upmconstant$ which then remains constant. 
	In this section, we provide intuition as to the meaning of the parameter $\upmconstant$ and show the effect of $\upmconstant$ on the convergence of the \acrshort{upm}.
	
	\subsubsection{Theoretical Results}
	\label{sec:supm_theory}

	
	In the previous section, we defined $\bracketupdate{\supmtag}$, the update function for the \acrshort{supm}, and showed in \Cref{lem:upm_maintains_bracket} that it maps an extended bracket which satisfies \Cref{eq:x_condition_1,eq:x_condition_2} to a new extended bracket which satisfies the same conditions by inserting the point $\bracketstep{\supmtag}(\bracket)$.
	However, an assumption made in the statement of \Cref{lem:upm_maintains_bracket} was that $\bracketstep{\supmtag}(\bracket) \in (\extbracketleftone,\extbracketmiddle)\cup(\extbracketmiddle, \extbracketrightone)$ .  
	We now show that there exists a constant $\upmconstant$ and function  $\upmscaling(\bracket)$ such that $\bracketstep{\supmtag}(\bracket) \in (\extbracketleftone,\extbracketmiddle)\cup(\extbracketmiddle, \extbracketrightone)$ holds.
	First we need the following preliminary results.
\fi


\begin{definition}
	A function $f:[a,b]\rightarrow \mathbb{R}$ is called unimodal if there exists $\zeta \in (a,b)$ such that $f$ is monotonically  decreasing on $(a, \zeta)$ and monotonically increasing on $(\zeta, b)$. If a function is not unimodal, then it is multi-modal.
\end{definition}

\begin{definition}
\label{def:nice_piecewise_smooth}
Let $f$ be a piece-wise smooth function. We call $f$ \textbf{\nice} if for any local minimum $x^*$ of $f$, there exists an open set $(a,b)\ni x^*$ such that $f$ is unimodal on $(a,b)$.
\end{definition}

When constructing convergence results, there are three levels of assumptions we might make. 
First is the case where $f$ is a piece-wise smooth function, $\bracket$ is a bracket, and nothing more is assumed. 
In this general context, we can do nothing beyond defining a minimum distance between points $\delta$.
If $|\bracketstep{\supmtag}(\bracket) - \extbracketmiddle| < \mindist$, we replace $\bracketstep{\supmtag}(\bracket)$ with $\argmin_x\; |x-\bracketstep{\supmtag}(\bracket)|$ such that $|x-\extbracketmiddle|\geq \mindist$, $\extbracketrightone-x\geq \mindist$ and $x-\extbracketleftone\geq \mindist$.
If $\bracketstep{\supmtag}(\bracket)=\extbracketmiddle$, then there exists two equally valid solutions: $\extbracketmiddle+\delta$ and $\extbracketmiddle-\delta$. 
When this applies, we arbitrarily set $\bracketstep{\supmtag}(\bracket)=\extbracketmiddle-\delta$.
By doing this, we are sure to converge eventually, even if slowly. 

For the next level of assumptions, we additionally require  $f$ to be in the form $f(x) = \max(f_L(x), f_R(x))$, and $\bracket$ satisfy $f(\extbracketleftj) = f_L(\extbracketleftj),\;j=1,2,3$ and $f(\extbracketrightj) = f_L(\extbracketrightj),\;j=1,2,3$.
Given these assumptions, we show in \Cref{lem:alpha_threshold} and \Cref{cor:supm_interior_selection} that $\bracketstep{\supmtag}(\bracket; \upmconstant) \in (\extbracketleftone, \extbracketrightone)$ given a sufficiently large choice of $\upmconstant$. 
In short, when there are not too many kinks in one place, then the \acrshort{upm} selects sensible points to evaluate $f$ at.

Finally, we add the assumptions that $f$ is unimodal and $\bracketouterdiam{\bracket}$ is sufficiently small. 
Under these circumstances, we show in \Cref{lem:supm_step_position} and \Cref{thm:supm_convergence_ns_quality2} that the distance between $\bracketstep{\supmtag}(\bracket)$ and the true solution $\truesol$ can be bounded.

Taken in the context of \nice functions, these results will imply that the \acrshort{supm} is stable when the function is not unimodal, and fast when it is. 
This is sufficient for our purposes given that a stable algorithm applied to a \nice function will eventually converge to a bracket in which the function is unimodal.



For the analysis that follows, we require the following external result:
%
\begin{theorem}[\cite{brentalgorithms}]
	\label{thm:divided_difference_accuracy}
	Suppose that $k, n\geq0$; $f \in C^{n+k}[a,b]$; $\zeta \in [a,b]$; and $x_0, \ldots, x_n$ are distinct points in $[a,b]$. Then
	\begin{equation}
	\label{eq:divided_difference_accuracy}
	\begin{aligned}
	f[x_0, \ldots, x_n] &= \frac{f^{(n)}(\zeta)}{n!}+\left(\sum \limits_{0\leq r_1\leq n} (x_{r_1}-\zeta) \right)\frac{f^{(n+1)}(\zeta)}{(n+1)!}\\
	&+\ldots+ \left( \sum \limits_{0 \leq r_1 \leq \ldots \leq r_k \leq n} \prod \limits_{j =1}^k (x_{r_j} - \zeta)\right) \frac{f^{(n+k)}(\zeta)}{(n+k)!} + E
	\end{aligned}
	\end{equation}
	where
	\begin{equation}
	\label{eq:divided_difference_error}
	E = \frac{1}{(n+k)!} \left( \sum \limits_{0 \leq r_1 \leq \ldots \leq r_k \leq n} \left(\prod \limits_{j =1}^k (x_{r_j} - \zeta)\right)[f^{(n+k)}(\xi_{r_1, \ldots, r_k}) - f^{(n+k)}(\zeta)]\right),
	\end{equation}
	and $\xi_{r_1, \ldots, r_k}$ are points in the interval spanned by $x_{r_1}, \ldots, x_{r_k}$ and $\zeta$.
\end{theorem}


\begin{lemma}
	\label{lem:alpha_threshold}
	Let $f$ be a function of the form $f = \max (f_L, f_R)$ such that $f_L''$ and $f_R''$ are both Lipschitz continuous with Lipschitz constants $M_L$ and $M_R$ respectively, and $\bracket$ be a an extended bracket of $f$  where $f(x^k_j) = f_k(x^k_j)$ for $j=1,2,3$ and $k \in \{L,R\}$. 
	If $ \upmconstant \geq \tfrac{1}{2}max(M_L, M_R)$ and $\upmscaling(\bracket) = \max(\extbracketrightthree - \extbracketleftone, \extbracketrightone - \extbracketleftthree)$, then $q^R(x;\bracket, \alpha)< f_R(x)\;\forall x \in [\extbracketleftone, \extbracketrightone)$ and $q^L(x;\bracket, \alpha)< f_L(x)\; \forall x \in (\extbracketleftone, \extbracketrightone]$ .
\end{lemma}
\begin{proof}
	Suppose $y \in (\extbracketleftone, \extbracketrightone]$.
	Since $f(\extbracketleftj) = f_L(\extbracketleftj)$ for $j=1,2,3$, we have:
	\begin{align*}
	&f_L(y) > q^L(y; \bracket, \upmconstant),\\
	\Leftrightarrow	&f_L(y) > f_L(\extbracketleftone) + f_L[\extbracketleftone, \extbracketlefttwo](y - \extbracketleftone) + (f_L[\extbracketleftone, \extbracketlefttwo, \extbracketleftthree]-\upmconstant \upmscaling(\bracket))(y - \extbracketleftone)(y - \extbracketlefttwo),\\
	\Leftrightarrow& f_L[y, \extbracketleftone] > f_L[\extbracketleftone, \extbracketlefttwo] + (f_L[\extbracketleftone, \extbracketlefttwo, \extbracketleftthree]-\upmconstant \upmscaling(\bracket))(y - \extbracketlefttwo),\\
	\Leftrightarrow &f_L[y, \extbracketleftone, \extbracketlefttwo] > f_L[\extbracketleftone, \extbracketlefttwo, \extbracketleftthree] - \upmconstant \upmscaling(\bracket),\\
	\Leftrightarrow& \upmconstant \upmscaling(\bracket) > f_L[\extbracketleftone, \extbracketlefttwo, \extbracketleftthree] - f_L[y, \extbracketleftone, \extbracketlefttwo].
	\end{align*}
	Analogously, for $y \in [\extbracketleftone, \extbracketrightone)$ we have:
	\[\upmconstant \upmscaling(\bracket) > f_R[\extbracketrightone, \extbracketrighttwo, \extbracketrightthree] - f_R[y, \extbracketrightone, \extbracketrighttwo].\]
	
	From \Cref{thm:divided_difference_accuracy}, we know that there exists $\xi^R_1\in [\extbracketrightone, \extbracketrightthree]$ such that $\tfrac{1}{2}f''(\xi^R_1) = f_R[\extbracketrightone, \extbracketrighttwo, \extbracketrightthree]$ and $\xi^R_2 \in [y, \extbracketrighttwo]$ such that $\tfrac{1}{2}f''(\xi^R_2) = f_R[y, \extbracketrightone, \extbracketrighttwo]$.
	Equivalently, the points $\xi^L_1$ and $\xi^L_2$ exist for $f_L[\extbracketleftone, \extbracketlefttwo, \extbracketleftthree]$ and $f_L[y, \extbracketleftone, \extbracketlefttwo]$ respectively. 
	Therefore:
	\begin{align*}
	 f_L[\extbracketleftone, \extbracketlefttwo, \extbracketleftthree] - f_L[y, \extbracketleftone, \extbracketlefttwo] &= \tfrac{1}{2} \left(f_L''(\xi^L_1) - f_L''(\xi^L_2) \right),\\
	 &\leq \tfrac{1}{2}M_L |\xi^L_1 - \xi^L_2|,\\
	 &\leq \tfrac{1}{2}M_L|\max(y, \extbracketleftone) - \extbracketleftthree |.
	\end{align*}
	The bound on the RHS depends on the value of $y$, and is itself bounded by $\tfrac{1}{2}M_L |\extbracketrightone - \extbracketleftthree|$ when considering $y \in [\extbracketleftone, \extbracketrightthree]$.
	Therefore, a sufficient condition for $f_L(y) > q^L(y; \bracket, \upmconstant)\; \forall y \in (\extbracketleftone, \extbracketrightone]$ is that $\upmconstant \upmscaling(\bracket) > \tfrac{1}{2}M_L |\extbracketrightone - \extbracketleftthree|$.  
	Similarly, a sufficient condition for $f_R(y) > q^R(y; \bracket, \upmconstant)\; \forall y \in [\extbracketleftone, \extbracketrightone)$ is $\upmconstant \upmscaling(\bracket) > \tfrac{1}{2}M_R |\extbracketrightthree - \extbracketleftone|.$
	
	Choosing $\upmscaling(\bracket) = \max(\extbracketrightone - \extbracketleftthree, \extbracketrightthree - \extbracketleftone)$, and $2 \alpha > \max(M_L, M_R)$ ensures that both of these are satisfied.
\end{proof}

\begin{corollary}
	\label{cor:supm_interior_selection}
	Let $f$ be a function of the form $f = \max (f_L, f_R)$ such that $f_L''$ and $f_R''$ are both Lipschitz continuous with Lipschitz constants $M_L$ and $M_R$ respectively, and $\bracket$ be a an extended bracket of $f$  where $f(x^k_j) = f_k(x^k_j)$ for $j=1,2,3$ and $k \in \{L,R\}$. 
	If $\upmconstant \geq \tfrac{1}{2} \max(M_L, M_R)$ and $\upmscaling(\bracket) = \max(\extbracketrightone - \extbracketleftthree, \extbracketrightthree - \extbracketleftone)$, then  $\bracketstep{\supmtag}(\bracket) \in (\extbracketleftone, \extbracketrightone)$. 
\end{corollary}
\begin{proof}
	As \Cref{lem:alpha_threshold} is applicable, we know that $q^R(x)< f_R(x) \;\forall x \in[\extbracketleftone, \extbracketrightone)$ and $q^L(x)<f_L(x)\;\forall x \in (\extbracketleftone, \extbracketrightone]$.
	Since $f = \max(f_L, f_R)$, it follows that $q^L(\extbracketmiddle) < f(\extbracketmiddle)$ and $q_R(\extbracketmiddle) < f(\extbracketmiddle) $. 
	For convenience, write $q(x) = \max(q^L(x;\bracket, \upmconstant), q^R(x;\bracket, \upmconstant))$ . Starting from \Cref{eq:UPM_gp}, we obtain: 
	\[ q(\bracketstep{\supmtag}(\bracket;\upmconstant)) \leq   q(\extbracketmiddle) < f(\extbracketmiddle) \leq \min(f(\extbracketleftone), f(\extbracketrightone)) = \min(q^L(\extbracketleftone), q^R(\extbracketrightone)).  \]
	But $q(\extbracketleftone) = q^L(\extbracketleftone)$ and $q(\extbracketrightone) = q^R(\extbracketrightone) $. 
	Therefore $ q(\bracketstep{\supmtag}(\bracket;\upmconstant))<\min(q(\extbracketleftone), q(\extbracketrightone))$, which implies that $ \bracketstep{\supmtag}(\bracket;\upmconstant) \neq \extbracketleftone$ or $\extbracketrightone$. 
	Since $\bracketstep{\supmtag}(\bracket;\upmconstant) \in [\extbracketleftone, \extbracketrightone]$ by definition (see \Cref{eq:UPM_gp}), it follows that $\bracketstep{\supmtag}(\bracket;\upmconstant) \in (\extbracketleftone, \extbracketrightone)$.
%
%
\end{proof}
\begin{remark}
	Between \Cref{cor:supm_interior_selection} and the requirement that the minimum distance between any two points in $\bracket$ is $\delta$, we are assured that the conditions for \Cref{lem:upm_maintains_bracket} will be satisfied.
\end{remark}

For the remainder of this section, we work towards bounding the convergence of the \acrshort{supm} when $f$ is piecewise-smooth and unimodal.
We define $\truesol$ to be the unique local minimum of $f$ in $(\extbracketleftone, \extbracketrightone)$.
When assessing the \acrshort{supm}'s convergence rate, we  consider two cases based on the values of $f'_L(\truesol) $ and $ f'_R(\truesol)$.

If $f'_R(\truesol) =0$ or $ f'_L(\truesol)=0$ holds, then we expect the \acrshort{supm} to behave similarly to Brent's Method. 
For example, if $f'_R(\truesol) =0$ and $f'_L(\truesol)<0$, then we expect $q^R$ to behave similarly to the interpolated polynomial from Brent's Method. 
The \acrshort{supm} will be slower however for two reasons. 
\begin{enumerate}
	\item When the \acrshort{supm} computes a new point, that value may be stored in $\extbracketmiddle$ which does not impact the next computation but only the one after that. 
	\item As only $q^R$ is converging, any point which is used to construct $q^L$ will play little to no role in determining $\bracketstep{\supmtag}(\bracket)$, and is therefore useless.
\end{enumerate}

We focus our attention primarily on the case where $f'_L(\truesol)<0$ and $f'_R(\truesol)>0$. 
In order to bound $|\bracketstep{\supmtag}(\bracket)|$, we need to understand what type of point $\bracketstep{\supmtag}(\bracket)$ returns. 
We address this in \Cref{lem:supm_step_position}.
First however, we introduce a shorthand notation for Divided Differences which we use for the remainder of this section.
\begin{align*}
f^k_1&= f(x^k_1),\; k \in \{L, R\},\\
f^k_2&= f[x^k_1, x^k_2],\; k \in \{L,R\},\\
f^k_3 &= f[x^k_1, x^k_2, x^k_3],\; k \in \{L,R\}.
\end{align*}

\begin{lemma}
	\label{lem:supm_step_position}
	Let $f$ be a piecewise smooth unimodal function of the form $f = \max (f_L, f_R)$ such that $f_L''$ and $f_R''$ are both Lipschitz continuous with Lipschitz constants $M_L$ and $M_R$ respectively, and $\bracket$ be a an extended bracket of $f$. 
	Further let $\upmconstant$ be given such that $2 \upmconstant \geq \max(M_L, M_R)$.
	If $f_L'(0)<0$ and $f_R'(0)>0$, then there exists $\epsilon>0$ such that if $\bracket \in (-\epsilon, \epsilon)^7$,
	then the model functions $q^L$ and $q^R$ intersect exactly once within $(\extbracketleftone, \extbracketrightone)$ and $\bracketstep{\supmtag}(\bracket)$ returns this unique intersection point. 
\end{lemma}
\begin{proof}
	Since $f= \max(f_L, f_R)$ is unimodal, and $\bracket$ is an extended bracket, it follows that $f(x^k_j) = f_k(x^k_j)$ for $j=1,2,3$ and $k \in \{L,R\}$. 
	We use this along with \Cref{lem:alpha_threshold} to show: $q^L(\extbracketleftone) = f(\extbracketleftone) = f_L(\extbracketleftone) \geq f_R(\extbracketleftone) >q^R(\extbracketleftone)$, and similarly $q^R(\extbracketrightone)> q^L(\extbracketleftone)$.
	By the Intermediate Value Theorem \cite[Theorem 4.23]{rudin1976principles}, it follows that $q^L$ intersects with $q^R$ at least once in $(\extbracketleftone, \extbracketrightone)$. 
	Moreover, since $q^L$ and $q^R$ are both quadratic functions, they must intersect exactly once.
	
	A sufficient but not necessary condition for the result to follow would be $\exists \epsilon>0$ such that $q^L$ and $q^R$ are monotonically decreasing and increasing respectively. 
	Once this holds, $\bracketstep{\supmtag}(\bracket) = \argmin_x \max(q^L(x), q^R(x))$ must refer to the intersection of $q^L$ and $q^R$.
	
	To show that $q^L$ is monotonically decreasing, it is sufficient to prove that $(q^L)'(\extbracketleftone)\leq 0$ and $(q^L)'(\extbracketrightone)< 0$. 
	Once this is shown, we are done since $(q^L)'$ is a linear function of $x$.	
	We find :
	\[ (q^L)'(\extbracketleftone) = f^L_2 = \frac{f(\extbracketleftone) - f(\extbracketlefttwo)}{\extbracketleftone - \extbracketlefttwo} \leq 0, \]
	since $f$ is unimodal and $\bracket$ is an extended bracket of $f$ (recall \Cref{def:extended_bracket}).
	For $(q^L)'(\extbracketrightone)$  we begin with 
	\begin{equation}
	\label{eq:minor_equation1}
	(q^L)'(\extbracketrightone)  = f^L_2 + (f^L_3 - \upmconstant \upmscaling(\bracket))(2 \extbracketrightone - \extbracketleftone - \extbracketlefttwo).
	\end{equation}   
	Next we apply \Cref{thm:divided_difference_accuracy} to conclude:
	\begin{align*}
	\exists \xi^L_1, \xi^L_2 \in [\extbracketlefttwo, \extbracketleftone] \;\text{such that}\;& f^L_2 = f_L'(0) + \tfrac{1}{2}(\extbracketleftone  f_L''(\xi^L_1)+ \extbracketlefttwo f_L''(\xi^L_2)),\\
	\exists \xi^L_3 \in [\extbracketleftthree, \extbracketleftone]\;\text{such that}\;& f^L_3 = \tfrac{1}{2}f''_L(\xi^L_3).
	\end{align*}
	From this and \Cref{eq:minor_equation1} it follows that
	\begin{align*}
	(q^L)'(\extbracketrightone)  &= f_L'(0) + \tfrac{1}{2} \left(  \extbracketleftone f_L''(\xi^L_1) + \extbracketlefttwo f_L''(\xi^L_2)\right) \\
	& + \tfrac{1}{2}(2\extbracketrightone - \extbracketleftone - \extbracketlefttwo)\left(f_L''(\xi^L_3) - 2\upmconstant \upmscaling(\bracket) \right).
	\end{align*}
	Since $f_L''$ and $f_R''$ are Lipschitz continuous, we know that they are bounded on $[\extbracketleftthree, \extbracketrightthree]$. 
	Let $M$ be constant which bounds $f_L''$ and $f_R''$. 
	We use this along with the fact that $\bracket \in (-\epsilon, \epsilon)^7$ to conclude:
	\begin{align*}
	(q^L)'(\extbracketrightone) &\leq f_L'(0) + \tfrac{1}{2} \left( - \extbracketleftone M - \extbracketlefttwo M + (2\extbracketrightone - \extbracketleftone - \extbracketlefttwo)(M - 2 \upmconstant \upmscaling(\bracket) \right),\\
	(q^L)'(\extbracketrightone)  & \leq f_L'(0) + \tfrac{1}{2} \left( 2 \epsilon M + 4 \epsilon (M - 2\upmconstant \upmscaling(\bracket) \right),\\
	& \leq f_L'(0) +\epsilon \left( 3M  - 4 \upmconstant \upmscaling(\bracket) \right) \leq f'_L(0) + 3 \epsilon M,
	\end{align*}
	since $0<\upmscaling(\bracket) <2\epsilon$ and $\upmconstant$ is constant.
	Therefore a sufficient condition for $(q^L)'(\extbracketrightone)<0$ to hold is that $\epsilon < |f_L'(0)|/3M$. 
	That $q^R$ is monotonically increasing given sufficiently small $\epsilon$ follows from the equivalent argument.

\end{proof}


Having established the conditions under which we can express  $\bracketstep{\supmtag}(\bracket)$ analytically, we now bound the distance between this point and the true minimiser.

\begin{theorem}
	\label{thm:supm_convergence_ns_quality2}
	Let $f:[a,b]\rightarrow \mathbb{R}$ be a piecewise smooth unimodal function of the form $f = \max (f_L, f_R)$ such that $\truesol\in [a,b]$ is a minimiser of $f$, $f_L'(\truesol)<0$, $f_R'(\truesol)>0$, and $\bracket$ be an extended bracket of $f$. If the functions $f_L''$ and $f_R''$ are both Lipschitz continuous (with constant $M_L$ and $M_R$) and $\upmconstant > \tfrac{1}{2}\max(M_L, M_R)$, then  there exists  $\epsilon>0$ such that
	\begin{align*} 
	|\bracketstep{\supmtag}(\bracket) | &\leq  \frac{\sqrt{2}}{f_R'(\truesol) - f_L'(\truesol)} \left( |\extbracketleftone-\truesol||\extbracketlefttwo-\truesol|(|\extbracketleftone-\truesol| + |\extbracketleftthree-\truesol|)M_L\right. \\
	&+ \left. |\extbracketrightone-\truesol||\extbracketrighttwo-\truesol|(|\extbracketrightone-\truesol| + |\extbracketrightthree-\truesol|)M_R \right), 
	\end{align*}
	when $\bracket \in  (-\epsilon, \epsilon)^7$.
\end{theorem}
\begin{proof}
	Since  \Cref{lem:supm_step_position} applies, we know that $\bracketstep{\supmtag}(\bracket)$ returns the unique point of intersection between $q^L$ and $q^R$ in $(\extbracketleftone, \extbracketrightone)$. 
	The intersection points of these quadratics are the roots to the polynomial $ax^2 + bx + c=0$, where:
	%
	%
	\begin{align*}
	a&=f^R_3 - f^L_3,\\
	b &= f^R_2 - f^L_2 - (\extbracketrightone +\extbracketrighttwo)(f^R_3-\upmconstant h) + (\extbracketleftone +\extbracketlefttwo)(f^L_3-\upmconstant h),\\
	c&=f^R_1 - f^L_1- \extbracketrightone f^R_2 + \extbracketleftone f^L_2 + \extbracketrightone \extbracketrighttwo (f^R_3-\upmconstant h) - \extbracketleftone \extbracketlefttwo (f^L_3-\upmconstant h).
	\end{align*}
	Assuming without loss of generality $\truesol = 0$, we invoke \Cref{thm:divided_difference_accuracy} to make the following substitutions for $\ell \in \{L, R\}$.
	\begin{align*}
	f^\ell_1 &= f_\ell(0) + x^\ell_1 f_\ell'(0) + \tfrac{1}{2}(x^\ell_1)^2 f_\ell''(\xi^\ell_1),&\xi^\ell_1 \in [0, x^\ell_1] , \\	
	f^\ell_2 &= f_\ell'(0) + \tfrac{1}{2}\left( x^\ell_1 f_\ell''(\xi^\ell_2) + x^\ell_2 f_\ell''(\xi^\ell_3) \right), & \xi^\ell_2 \in [0, x^\ell_1], \; \xi^\ell_3 \in [0, x^\ell_2],\\
	f^\ell_3 &= \tfrac{1}{2} f_\ell''(\xi^\ell_4), & \xi^\ell_4 \in [0, x^\ell_3].
	\end{align*}
	By inserting the substitutions above into the definitions of $a,b$ and $c$, we get:
	%
	\begin{align*}
	a&= \tfrac{1}{2}\left(f_R''(0) - f_L''(0) \right) + K^R_a - K^L_a,\\
	b&= f_R'(0) - f_L'(0) + K_b^R - K_b^L, \\
	c&= K^R_c - K^L_c,
	\end{align*}
	where $K^\ell_a$, $K^\ell_b$ and $K^\ell_c$ are in turn defined for $\ell \in \{L, R\}$ by:
	\begin{align*}
	K^\ell_a &= \tfrac{1}{2} \left(f_\ell''(\xi^\ell_4) - f_\ell''(0) \right), \\
	K^\ell_b &= \tfrac{1}{2}\left(x^\ell_1 (f_\ell''(\xi^\ell_2) - f_\ell''(\xi^\ell_4)) + x^\ell_2(f_\ell''(\xi^\ell_3) - f_\ell''(\xi^\ell_4)) \right) + (x^\ell_1 + x^\ell_2)\upmconstant \upmscaling(\bracket),\\
	K^\ell_c &= \tfrac{1}{2} x^\ell_1 \left( x^\ell_1 (f_\ell''(\xi^\ell_1) - f_\ell''(\xi^\ell_2) )+ x^\ell_2(f_\ell''(\xi^\ell_4) - f_\ell''(\xi^\ell_3)) \right) - x^\ell_1 x^\ell_2 \upmconstant \upmscaling(\bracket).
	\end{align*}
	Writing $a,b$ and $c$ in terms of these $K$ values allows us to give bounds using the assumption that $f_L''$ and $f_R''$ are Lipschitz Continuous . 
	\begin{subequations}
	\begin{align}
		|K^\ell_a| &\leq  \tfrac{1}{2}M_\ell |x^\ell_3| \leq \tfrac{1}{2}M_\ell \epsilon , \\
		|K^\ell_b| &\leq \tfrac{1}{2}(|x^\ell_1| + |x^\ell_2|)\left(M_\ell|x^\ell_3|   + \upmconstant \upmscaling(\bracket) \right) \leq \epsilon^2(M_\ell + 2 \upmconstant),\\
		|K^\ell_c| &\leq \tfrac{1}{2} |x^\ell_1||x^\ell_2| M_\ell \left( |x^\ell_1| +  |x^\ell_3| \right) \leq M_\ell \epsilon^3.\label{pfeq:K_bound}
		\end{align}
	\end{subequations}

	Choose an $\epsilon_1$ small enough to ensure $\tfrac{1}{2}(f_R'(0) - f_L'(0))> K^L_b - K^R_b$ (note that $b \rightarrow f_R'(0) - f_L'(0) $ as $\epsilon \rightarrow 0$). 
	If $\epsilon \leq \epsilon_1$, then it follows that $b>0$ which in turn implies:
	\begin{equation}
	\label{eq:analytic_step}
	\bracketstep{\supmtag}(\bracket) = \frac{-b + \sqrt{b^2 - 4ac}}{2a},
	\end{equation} 
	because the other root lies outside $[\extbracketleftone, \extbracketrightone]$ as $\epsilon \rightarrow 0$.
	We apply Taylors Theorem to conclude $\exists \xi \in [-|4a c|, |4a c|]$ such that $\sqrt{b^2 - 4ac}= b + 2ac/\sqrt{b^2 - \xi}$.
	Since $ac \rightarrow 0$ as $\epsilon \rightarrow 0$, there exists $\epsilon_2$ such that $|4ac| \leq \tfrac{1}{2}b^2$ when $\epsilon < \epsilon_2$. 
	Then if $\epsilon < \min(\epsilon_0, \epsilon_1, \epsilon_2)$, we have:
	\begin{align*}
	|\bracketstep{\supmtag}(\bracket)|& = \left| \frac{-b + \sqrt{ b^2 - 4ac } }{2 a}  \right| = \cfrac{1}{2|a|}\left|\cfrac{2ac}{\sqrt{b^2 - \xi} }\right|, & \xi \in [0, 4ac],\\
	&\leq \frac{|c|}{\sqrt{b^2 - \xi}}\leq \frac{|c|}{\sqrt{b^2 - 4|ac|}}\leq \sqrt{2}\left|\frac{c}{b}\right|, \\
	&= \sqrt{2} \frac{K^R_c - K^L_c}{f_R'(0) - f_L'(0) + K^R_b - K^L_b} \leq 2\sqrt{2} \frac{|K^R_c|+ |K^L_c|}{f_R'(0) - f_L'(0)}.
	\end{align*}
	Finally we insert \Cref{pfeq:K_bound} and the result follows.
\end{proof}

\Cref{thm:supm_convergence_ns_quality2} is remarkably similar to \cite[Theorem 4.1]{mifflin1990superlinear}, the equivalent result which bounds the convergence of \algmiff.
Using our own notation for their result, Mifflin and Strodiot showed $|\extbracketmiddle_i - x^*| \leq A_i |x^* - \extbracketleft{1}{\idxi}| + B_i |\extbracketright{1}{\idxi}- x^*|$, for each iteration $i$, where $A_i$ and $B_i$ are sequences which converge to $0$ as $i \rightarrow \infty$.
If we rewrite the result from \Cref{thm:supm_convergence_ns_quality2} in the same form, the sequence $A_i$ may be written explicitly in terms of $\extbracketleft{1}{\idxi}, \extbracketleft{2}{\idxi}, $ and $\extbracketleft{3}{\idxi}$, and similarly for $B_i$.
Then, if $\extbracketright{1}{\idxi}$ and $\extbracketleft{1}{\idxi}$ both converge to $x^*$, the rate at which $A_i$ and $B_i$ tend to zero begins to resemble quadratic convergence. 

However, \Cref{thm:supm_convergence_ns_quality2}, like \cite[Theorem 4.1]{mifflin1990superlinear} does not prove quadratic convergence, nor even linear convergence. 
This is because we cannot know a priori whether the sequences $\{\extbracketleft{1}{\idxi}\}$ and $\{\extbracketright{1}{\idxi}\}$ converge at all, nor how quickly if they do.
Suppose there exists an iteration $i_0$ such that $\forall i>i_0$, $\extbracketright{1}{\idxi}$ is constant.
In this case, the bound given by \Cref{thm:supm_convergence_ns_quality2} will decrease, but not converge to $0$. 
In practice, when this occurs $\extbracketmiddle_i$ will still converge to $x^*$ but only sub-linearly, while $\bracketinnerdiam{\bracket_i}\not \rightarrow 0$.

On the other hand, suppose that both $\{\extbracketleft{1}{\idxi}\}$ and $\{\extbracketright{1}{\idxi}\}$ converge to $x^*$, but we only update $\extbracketright{1}{\idxi}$ once every $k$ iterations.
Then we may say that $|\extbracketmiddle_{k i} - x^*|\rightarrow0$ super-linearly as $i \rightarrow \infty$, but depending on the size of $k$, this may not be particularly useful.

In short, there are three main weaknesses of the \acrshort{supm} which lead to its failure. 
The first is that we need to choose $\upmconstant$, a suitable value of which we cannot know in advance. 
The second is that convergence only occurs in reasonable time if both updates of $\extbracketright{1}{\idxi}$ and $\extbracketleft{1}{\idxi}$ occur regularly. 
Finally, we require a bracket sufficiently close to the solution that $f$ is unimodal before any of our results apply.
These shortcomings of the \acrshort{supm} motivate subsequent sections, in which we describe more stable variations of the \acrshort{upm}.

%
%


\ifnum \condensed=0

	\subsubsection{Numerical Results}
	\label{sec:supm_testing}
	Having discussed the theoretical properties of the \acrshort{supm}, in this section we now show it in action.
	We will demonstrate how its performance is affected by choice of $\upmconstant$ and on the objective function. 
	Before showing the results of our numerical experiment, we first provide specific examples to show what occurs in practice. 
	
	When testing the \acrshort{supm}, we apply it to three categories of test-functions: smooth unimodal, non-smooth unimodal and smooth multimodal. 
	To our knowledge there is not a commonly used set of univariate test functions. 
	Therefore, we have taken some test functions from Jamil et al \cite{jamil2013literature}, keeping the ones which are well defined for one variable. 
	Adding to this, we have added a few of our own design to be adversarial examples based on our observations. 
	We scale each test function $f$ such that the minimum Lipschitz constant for $f'''$ lies between $0.8$ and $1$.
	For the category of smooth unimodal, we use the following  test functions:
	\begin{align}
		f^{SU}_1 &= -\frac{1}{\sqrt{e}}\exp\left( -\tfrac{1}{2} x^2 \right),& -1 \leq x \leq 1,\\
		f^{SU}_2 &= \frac{1}{24} x^4 ,& -1 \leq x \leq 1,\\
		f^{SU}_3 &= \frac{1}{11}\left(- \sin (2x - \tfrac{1}{2} \pi) - 3 \cos x - \tfrac{1}{2} x \right),& -2.5 \leq x \leq 3,\\
		f^{SU}_4 &= \frac{1}{2500}\left( \frac{x^2}{2} -\frac{ \cos \left(5 \pi x\right) }{25 \pi^2} - \frac{x \sin(5 \pi x)}{5 \pi } \right), &-10\leq x \leq 10,\\ 
		f^{SU}_5 &= -\frac{1}{250} \left(x^{\frac{2}{3}} + (1-x^2)^{\frac{1}{3}}\right),&0.1 \leq x \leq 0.9,\\
		f^{SU}_6 &= \frac{1}{6000} \left( e^x + \frac{1}{\sqrt{x}}\right) , & 0.1 \leq x \leq 3,\\
		f^{SU}_7 &= -\frac{1}{13}(16 x^2 - 24 x + 5) e^{-x},&1.3 \leq x \leq 3.9.
	\end{align}
	\begin{figure}[h!]
		\centering
		\includegraphics[width=1.0\textwidth]{PythonFigs/testcaseSU1_outcome.pdf}
		\caption{We demonstrate the \acrshort{supm} with various values of $\upmconstant$ when applied to $f^{SU}_1(x)$.}
		\label{fig:SU1_outcome}
	\end{figure}
	
	\begin{table}
		\centering
		\begin{filecontents*}{CSVs/SUtest1000.csv}
0, 0.01, 0.1, 1, 10, 100, EUPM, DUPM, Brent, GS, Mifflin
3.267802444291763031e-01, 3.659892122724482388e-01, 4.272625181999492039e-01, 4.912230975015475787e-01, 5.472157051746484591e-01, 5.894920425674853304e-01, 6.191561293810334821e-01, 4.800351615480049072e-01, 2.158506906833654038e-01, 6.390302222532231458e-01, 1.476810240863207502e-01
6.800002516316703272e-01, 6.409133074038888545e-01, 6.164458669356442932e-01, 6.182722693065396680e-01, 6.187012688664940141e-01, 6.185920470852268593e-01, 6.187822190737763961e-01, 6.182067161922861764e-01, 3.281686223631667954e-01, 6.389634849389871363e-01, 6.613863929738110770e-01
inf, inf, 4.848267179735524790e-01, 5.437694844516764991e-01, 5.881877575592735496e-01, 6.186448077077280061e-01, 6.349266326649241066e-01, 5.292865976631807579e-01, 3.329082595298186176e-01, 6.375270801518500008e-01, inf
inf, 6.102362707646488138e-01, 6.104879186070979458e-01, 6.104381399577749612e-01, 6.103923070938810369e-01, 6.110293101280865891e-01, 6.106740718038189408e-01, 6.102681504112197342e-01, 5.625865860868042301e-01, 6.356752740661003598e-01, inf
inf, 4.342014792192867567e-01, 5.468523786461898517e-01, 5.925318873957589050e-01, 6.233064637576102296e-01, 6.327386040010936075e-01, 6.323170196868898030e-01, 4.853988318313961492e-01, 2.756170550952111720e-01, 6.407968065652881462e-01, 2.581089441677300589e-01
inf, 5.450237715117506321e-01, 5.870955909221332591e-01, 6.149958302397957599e-01, 6.246055990496257593e-01, 6.247458566140964287e-01, 6.248590855703404223e-01, 5.569075125264001391e-01, 2.812329331056069859e-01, 6.143217036111685569e-01, inf
3.500156492962426813e-01, 4.025580742655491817e-01, 4.965863092392512335e-01, 5.520505845584199323e-01, 5.944147960066219483e-01, 6.233377335781912221e-01, 6.352584772099363342e-01, 4.884537964543064770e-01, 2.574536644279001862e-01, 6.385692176713058110e-01, inf
\end{filecontents*}
\pgfplotstabletypeset[	col sep = comma,
		/pgf/number format/precision=4,
		create on use/Functions/.style={
			create col/set list={$f^{SU}_1$,$f^{SU}_2$,$f^{SU}_3$,$f^{SU}_4$,$f^{SU}_5$,$f^{SU}_6$,$f^{SU}_7$}
		},
		columns/Functions/.style={string type},
		columns={Functions,0,0.1,1,10,EUPM,Brent,Mifflin},
		display columns/1/.style={column name={$\alpha = 0$}},
		display columns/2/.style={column name={$\alpha = 0.1$}},
		display columns/3/.style={column name={$\alpha = 1$}},
		display columns/4/.style={column name={$\alpha = 10$}},
		display columns/5/.style={column name={EUPM}},
		display columns/6/.style={column name={Brent}},
		display columns/7/.style={column name={Mifflin}},
		every head row/.style={
			before row={\toprule}, 
			after row={\toprule}
		},
		every last row/.style={after row=\bottomrule}, 
		]{CSVs/SUtest1000.csv}
		\caption{Average convergence rate of different algorithms on the set of smooth unimodal test functions.}
		\label{tab:SUperformance}
	\end{table}
	The results obtained from running the \acrshort{supm} on $f^{S1}_6(x)$ are shown in \Cref{fig:SU6_outcome}.
	We see that Brent's method converges fastest, which is not surpising given that Brent's method is designed for such functions. 
	Next we see the \acrshort{supm} converging slightly slower than Brent's method for $\upmconstant=0,0.1,1$. 
	Of these, the means of convergence the choice of $\upmconstant=0$ is peculiar. 
	In particular, we see a period in which no progress is made to collapse the bracket, followed by a step which completes the convergence process instantly.
	
	In practice, we associate this behaviour with only one of the model functions converging to approximate the true function. 
	Therefore, the test points $\bracketstep{\supmtag}(\bracket_i)$ approach the true solution $\zeta$ from one side only. 
	This results in only one of $\extbracketleft{1}{\idxi}$ and $\extbracketright{1}{\idxi}$ converging to $\zeta$, the other remaining constant.
	However, this behaviour is only possible when the model function is converging to a smooth convex function. 
	Although the algorithm works in this case, from experience we associate its success with luck, and do not take much more from it.
	
	Finally we note that high values of $\upmconstant$ result in the \acrshort{supm} converging similarly to Golden Section.
	
	The next class of test functions we consider is that of non-smooth unimodal:
	\begin{align}
		f^{NU}_1 &= -60000 \exp\left(-\frac{|x|}{50} \right),& -32 \leq x \leq 32,\\
		f^{NU}_2 &= \frac{1}{6} \max\left(\frac{1}{x+3}, \log(x) \right),& -2 \leq x \leq 10,\\
		f^{NU}_3 &=\frac{1}{24} \max \left(\frac{1}{x+3}, \frac{1}{(x-3)^2} \right),& -2 \leq x \leq 2,\\
		f^{NU}_4 &= \frac{1}{160}\max\left(\frac{1}{x+3}, \exp(x) \right), &-2\leq x \leq 5,\\ 
		f^{NU}_5 &= \frac{1}{150}\max\left( \exp(-x), \exp(x) \right),&-5 \leq x \leq 5.
	\end{align}

	
	%
	\begin{figure}[h!]
		\centering
		\includegraphics[width=1.0\textwidth]{PythonFigs/testcaseNU2_outcome.pdf}
		\caption{We demonstrate the \acrshort{supm} with various values of $\upmconstant$ on the $f^{NU}_2$.}
		\label{fig:NU2_outcome}
	\end{figure}

	\begin{table}
		\centering
		\begin{filecontents*}{CSVs/NUtest1000.csv}
0, 0.01, 0.1, 1, 10, 100, EUPM, DUPM, Brent, GS, Mifflin
inf, inf, 1.868305444749769340e-01, 2.684311096107854233e-01, 3.272306863043013547e-01, 3.708289200135256314e-01, 6.188417717930393414e-01, 2.640225560859142395e-01, 5.115005356445111451e-01, 6.345140346172843948e-01, inf
inf, inf, 3.998837803039220207e-01, 4.445379506525618529e-01, 4.806353095852366475e-01, 5.120590952867639656e-01, 6.265373093470425481e-01, 4.270354356818293895e-01, 6.337059659464091554e-01, 6.364101032538217462e-01, inf
inf, inf, 4.537670428795743671e-01, 4.931247194250562016e-01, 5.254673611701987657e-01, 5.533917461545891925e-01, 6.413104660731403506e-01, 4.420990965557121677e-01, 6.457235445739026858e-01, 6.380097837701509400e-01, 3.259301021925230524e-01
inf, 4.212710370705753937e-01, 4.152834270366949077e-01, 4.665434255422860055e-01, 5.005704284943386373e-01, 5.302154729582594372e-01, 6.203998938230267379e-01, 4.050684044272888640e-01, 5.934306387578475173e-01, 6.253396755895377357e-01, 3.377983011049813244e-01
inf, 4.171920404856992448e-01, 4.203987182440860337e-01, 4.645000002791097593e-01, 4.974359132201823797e-01, 5.271951515881710959e-01, 6.194775858887014985e-01, 4.141885950657484550e-01, 4.903091623470886562e-01, 6.367320612107140176e-01, 3.592294375720715971e-01
\end{filecontents*}
\pgfplotstabletypeset[	col sep = comma,
		/pgf/number format/precision=4,
		create on use/Functions/.style={
			create col/set list={$f^{NU}_1$,$f^{NU}_2$,$f^{NU}_3$,$f^{NU}_4$,$f^{NU}_5$}
		},
		columns/Functions/.style={string type},
		columns={Functions,0,0.1,1,10,EUPM,Brent,Mifflin},
		display columns/1/.style={column name={$\alpha = 0$}},
		display columns/2/.style={column name={$\alpha = 0.1$}},
		display columns/3/.style={column name={$\alpha = 1$}},
		display columns/4/.style={column name={$\alpha = 10$}},
		display columns/5/.style={column name={EUPM}},
		display columns/6/.style={column name={Brent}},
		display columns/7/.style={column name={Mifflin}},
		every head row/.style={
			before row={\toprule}, 
			after row={\toprule}
		},
		every last row/.style={after row=\bottomrule}, 
		]{CSVs/NUtest1000.csv}
		\caption{Average convergence rate of different algorithms on the set of non-smooth unimodal test functions.}
		\label{tab:NUperformance}
	\end{table}
	We show the results of applying the \acrshort{supm} to $f^{NU}_2$ in \Cref{fig:NU2_outcome}. 
	Note that the \acrshort{supm}  fails to converge when $\upmconstant = 0$. 
	Excluding $\upmconstant=0$, we see that the larger $\upmconstant$ is, the slower the \acrshort{supm} converges. 
	While convergence occurs at different rates, the \acrshort{supm} appears to converge super-linearly for all finite values of $\upmconstant$. 
	Meanwhile, the performance of Brent's method is indistinguishable from Golden section, which in turn is marginally better than the limiting case where we take $\upmconstant \rightarrow \infty$. 
	We focus more on this final case in \Cref{sec:eupm}.
	
	Finally we list our smooth multimodal test functions:
	\begin{align}
		f^{SM}_1 &= \frac{x^6}{300} \left(2 + \sin\left(\frac{1}{x} \right)\right),& -1 \leq x \leq 1,\\
		f^{SM}_2 &= -\frac{1}{40000}\sin\left(5 \pi x \right)^6 ,& -1 \leq x \leq 1,\\
		f^{SM}_3 &= -\frac{1}{250000} \sin\left(5 \pi (x^{\frac{3}{4}}-\frac{1}{20}) \right)^6,& 0.01 \leq x \leq 1,\\
		f^{SM}_4 &= \frac{1}{5}\left( \sin\left(\frac{16}{15} x - 1\right)+ \sin \left(\frac{16}{15} x - 1\right)^2 \right),  &-1\leq x \leq 1,\\ 
		f^{SM}_5 &= \frac{x^2}{4000} -\cos x + 1,&100 \leq x \leq 100,\\
		f^{SM}_6 &=\frac{1}{71} \left( \left(\log (x-2) \right)^2 +\left(\log (10-x) \right)^2 - x^{\frac{1}{5}}\right), & 2.5 \leq x \leq 9.5,\\
		f^{SM}_7&= \frac{1}{40}  \left(\sin(x) + \sin \left(\frac{10 x}{3}\right) + \log(x) + \frac{21}{25} x \right), & 0.5\leq x \leq 10.
	\end{align}

	\begin{figure}[h!]
		\centering
		\includegraphics[width=1.0\textwidth]{PythonFigs/testcaseSM5_outcome.pdf}
		\caption{We demonstrate the \acrshort{supm} with various values of $\upmconstant$ on the function $f^{SM}_5$.  }
		\label{fig:SM5_outcome}
	\end{figure}

	\begin{table}
		\centering
		\begin{filecontents*}{CSVs/SMtest1000.csv}
0, 0.01, 0.1, 1, 10, 100, EUPM, DUPM, Brent, GS, Mifflin
inf, 6.417447209213930082e-01, 6.145538294522997491e-01, 6.137931730763693805e-01, 6.143601505316452771e-01, 6.136770479238801146e-01, 6.149399847048862000e-01, 6.223964193957512991e-01, 4.827908073819408563e-01, 6.391287237155460765e-01, inf
inf, inf, 5.076473248404567151e-01, 5.579151440922461846e-01, 5.960889957445194565e-01, 6.167074251245561189e-01, 6.182594029578646344e-01, 4.860131303907059275e-01, 3.587689684990479355e-01, 6.212019769184603524e-01, inf
inf, 4.846884638833482972e-01, 5.481325330788927586e-01, 5.918103288178249155e-01, 6.177394953272689060e-01, 6.233625018676751672e-01, 6.234205086822433017e-01, 5.085084300548894376e-01, 3.372016264353951964e-01, 6.187535245539960149e-01, inf
inf, inf, 4.535959404875268275e-01, 5.308438641935895319e-01, 5.789179915422920653e-01, 6.140703122379308487e-01, 6.337574015852992515e-01, 5.282278687514335713e-01, 2.703102859804504821e-01, 6.391180945675777325e-01, inf
inf, inf, 4.703796680941505493e-01, 5.187810971892090617e-01, 5.582084410665966168e-01, 5.893225946038999075e-01, 6.164562268469541140e-01, 5.142568100588714719e-01, 3.109510815490837965e-01, 6.278569443467145739e-01, inf
inf, 5.263376281885575603e-01, 5.598522405467014629e-01, 5.961887006976468451e-01, 6.190535028952102170e-01, 6.245048175538447932e-01, 6.240059726460778222e-01, 5.352530830861477185e-01, 2.982370183993352519e-01, 6.147915692839281965e-01, inf
inf, inf, 4.732989992150125613e-01, 5.353029776657690642e-01, 5.775549489213254173e-01, 6.085839109774059397e-01, 6.223935999543626085e-01, 4.814522428941768184e-01, 2.730028879916923668e-01, 6.114125297454715557e-01, inf
\end{filecontents*}
\pgfplotstabletypeset[	col sep = comma,
		/pgf/number format/precision=4,
		create on use/Functions/.style={
			create col/set list={$f^{SM}_1$,$f^{SM}_2$,$f^{SM}_3$,$f^{SM}_4$,$f^{SM}_5$,$f^{SM}_6$,$f^{SM}_7$}
		},
		columns/Functions/.style={string type},
		columns={Functions,0,0.1,1,10,EUPM,Brent,Mifflin},
		display columns/1/.style={column name={$\alpha = 0$}},
		display columns/2/.style={column name={$\alpha = 0.1$}},
		display columns/3/.style={column name={$\alpha = 1$}},
		display columns/4/.style={column name={$\alpha = 10$}},
		display columns/5/.style={column name={EUPM}},
		display columns/6/.style={column name={Brent}},
		display columns/7/.style={column name={Mifflin}},
		every head row/.style={
			before row={\toprule}, 
			after row={\toprule}
		},
		every last row/.style={after row=\bottomrule}, 
		]{CSVs/SMtest1000.csv}
		\caption{Average convergence rate of different algorithms on the set of smooth multimodal test functions.}
		\label{tab:SMperformance}
	\end{table}
	
	As in the example with the smooth unimodal test function, Brent's method is clearly the fastest. 
	Similarly to the non-smooth example, we see the \acrshort{supm} getting slower as $\upmconstant$ increases. 
	This is of course with the exception of $\upmconstant=0$ where the \acrshort{supm} fails to converge. 
	It is worth noting that even though none of the theory for the \acrshort{supm} was designed for multi modal function, it often converges anyway. 
	We discuss the possible reasons for this in \Cref{sec:eupm}.
	
	In the examples we have shown, there exist a few occasions for which the \acrshort{supm} has failed to converge.
	In order to show why this occurs, we present a specific example below, showing what properties $q^L$ and $q^R$ have when this happens.
	
	\begin{example}
		\label{ex:supm_failure}
		\begin{figure}
			\centering
			\includegraphics[width=1\textwidth]{PythonFigs/SUPM_Failure_6.pdf}
			\caption{When the \acrshort{supm} is run for 6 iterations on $f^{NU}_3$ from the starting point $\bracket_0$ (see \Cref{ex:supm_failure}), the model functions $q^L$ and $q^R$ are as plotted above .}
			\label{fig:supm_failure}
		\end{figure}
		Consider the objective function $f^{NU}_3$ with $\upmconstant = 0$ and
		\[\bracket_0 = (-1.811263, -2.171330, -2.23927, 1.820150, 2.102197, 2.293404, 2.334091).\]
		If we run the \acrshort{supm} up to the 6th iteration, we find the scenario illustrated in \Cref{fig:supm_failure}.

		In particular, we see a situation where the choice of $\tilde x$ is determined only by $q^L$, which in turn happens because $q^L(x)>q^R(x)$ for every $x \in [\extbracketleftone, \extbracketrightone]$. 
		The test point $\tilde x$ is the same as $\extbracketrightone$. 
		This means that in the following iteration, $q^R$ will be updated but not $q^L$.
		Thus the same thing will happen again, resulting the \acrshort{supm} stalling.
	\end{example}
	
	Our takeaway from this example is that the \acrshort{upm} process may fail if the model functions $q^L$ and $q^R$ do not underestimate the objective function $f$.
	This is a consequence of $\upmconstant$  not being large enough, and is consistent is \Cref{fig:NU4_outcome,fig:SM5_outcome}.

	Having chosen a set of test problems, all that remains is to run our various algorithms on each. 
	For each test function, we run Golden Section, Brent's Method and the \acrshort{supm} with $\upmconstant$ values $0.1, 1, 10, 100, 1000, \infty$ for randomly generated starting brackets.
	
	When constructing random brackets for test functions, it is important to note that we cannot guarantee that we will find a bracket which satisfies \Cref{eq:x_condition_2} for non-convex test functions. 
	Therefore, we relax this condition to:
	\begin{equation}
	\label{eq:bracket_condition_relaxed}
	 f(\extbracketleftone)\geq f(\extbracketmiddle)\leq f(\extbracketrightone).
	\end{equation}
	
	For each test function, we have an interval $[a,b]$ on which the function is intended to be used. 
	We sample 4 points uniformly from both $[a, a+\tfrac{1}{5}(b-a)]$ and $[b - \tfrac{1}{2}(b-a), b]$.
	After evaluating the function at these 8 points, we select the point where $f$ is minimised to be $\extbracketmiddle$; usually this will be one of the points closer to $\tfrac{1}{2}(b-a)$. 
	Finally we order the remaining points, and select the 3 closest on the left and right to be $\extbracketleftthree, \extbracketlefttwo, \ldots, \extbracketrightthree$. 
	
	Given that we have picked $\extbracketmiddle$ to be the point where $f$ is minimised, \Cref{eq:bracket_condition_relaxed} is satisfied automatically. 
	Furthermore, \Cref{eq:x_condition_1} is satisfied because of how we ordered the points.

	Having generated the bracket, we run each algorithm and compute the average convergence rate over 1000 different starting brackets.
	If the algorithm converged for a particular starting bracket, then its average convergence rate is by definition between 0 and 1, the smaller the better.
	If an algorithm did not converge, we assign $\infty$ instead.
	We do this because the purpose of the \acrshort{supm} is to be a robust algorithm. 
	Therefore we are not interested in it having a probability of success. 
	Either it always converges or it does not.
	We present the outcome of this experiment in \Cref{tab:supm_results}.
	
	From \Cref{tab:supm_results}, we see justification for our observations made earlier. 
	In particular, when the \acrshort{supm} converges, it usually converges faster for small values of $\upmconstant$.
	However, the smaller the value of $\upmconstant$, the less likely the \acrshort{supm} will converge in the first place. 
	Two notable exceptions to this are functions $f^{SU}_1$ and $f^{SM}_1$.
	A feature that these functions have in common is that the higher order derivatives go to zero near to the solution.

	For the smooth test functions, we see that Brent's Method is always the best. 
	However, for non-smooth functions it performs little better then Golden Section.
	On the other hand, the \acrshort{supm} performs between Brent's method and Golden Section for smooth functions and is far superior for non-smooth functions.
	Most interestingly, we see that \acrshort{supm} with $\upmconstant = \infty$ (we later name this the \acrlong*{eupm}) is effectively equivalent to Golden Section for all our experiments.
	

\fi	

\ifnum \condensed=0
	\subsection{Extremal Underestimating Polynomial Method}
\else
	\section{Extremal Underestimating Polynomial Method}
\fi
\label{sec:eupm}

In \ifnum \condensed=0 \Cref{sec:supm_theory} \else \Cref{sec:upm} \fi, we showed that if $\upmconstant$ was sufficiently large, then the $\bracketstep{\supmtag}$ is guaranteed to select a point within $(\extbracketleftone, \extbracketrightone)$.
Coupled with a restriction on the minimum distance between points, this is enough to ensure convergence eventually. 
However, for an arbitrary piecewise smooth objective function, we have no idea what a suitable $\upmconstant$ might be. 
As there is no finite value of $\upmconstant$ which satisfies \Cref{lem:alpha_threshold} for all functions, we ask what happens if we let $\upmconstant \rightarrow \infty$ and name the corresponding variation of the \acrshort{supm} the \acrfull*{eupm}, denoting its step function and update function by $\bracketstep{\eupmtag}$ and $\bracketupdate{\eupmtag}$ respectively.

%

As $\upmconstant \rightarrow \infty$ we are assured by \Cref{lem:alpha_threshold} that $q^L$ and $q^R$ will underestimate $f_L$ and $f_R$ respectively within $[\extbracketleftone, \extbracketrightone]$.
Furthermore, if $\upmscaling(\bracket)>0$, then we can pick a sufficiently large $\upmconstant$ such that the coefficient of $x^2$ in $q^k(x;\bracket,\alpha)$ (see \Cref{eq:UPM_model_function_definitions}) will be negative.
Once $\upmconstant$ satisfies both of these, then $q^L$ and $q^R$ will intersect at exactly one point, which is the point that  $\bracketstep{\supmtag}(\bracket)$ returns. 
We define $\bracketstep{\eupmtag}(\bracket)$ to be this point which we calculate by identifying  the root of $q^R(x; \bracket, \upmconstant) - q^L(x; \bracket, \upmconstant)=0$
, which remains bounded as $\upmconstant \rightarrow \infty$. 
By starting from \Cref{eq:analytic_step} and noting the Taylor Series as $\upmconstant \rightarrow \infty$, we find:

\begin{equation}
\label{eq:UPM_bigM_next_point}
\bracketstep{\eupmtag}(\bracket) = \frac{\extbracketrightone\extbracketrighttwo - \extbracketleftone \extbracketlefttwo}{\extbracketrightone+\extbracketrighttwo - \extbracketleftone - \extbracketlefttwo}.
\end{equation}

Note that the step function $\bracketstep{\eupmtag}$ does not depend on either $\extbracketleftthree$ or $\extbracketrightthree$. 
Therefore, within the context of this section, our extended bracket $\bracket$ is of the form:
\begin{equation}
\label{def:bracket_extremal}
\bracket = \left( \begin{array}{ccccc} \extbracketlefttwo&\extbracketleftone&\extbracketmiddle&\extbracketrightone&\extbracketrighttwo \end{array} \right)^T,
\end{equation} 

In order to construct $\bracketupdate{\eupmtag}$, we first define $\tilde{\bracketupdate{}}_{\supmtag}$ to be the update function $\bracketupdate{\supmtag}$, were we remove the entries corresponding to $\extbracketrightthree$ and $\extbracketleftthree$. 
Using this notation, we define:
\begin{equation}
\label{eq:eupm_update}
\bracketupdate{\eupmtag}(\bracket):= \tilde{\bracketupdate{}}_{\supmtag}(\bracketstep{\eupmtag}(\bracket), \bracket).
\end{equation}
Given $\bracketstep{\eupmtag}$ and $\bracketupdate{\eupmtag}$, the \acrshort{eupm} is precisely in the form of \Cref{alg:bracketing_method_template}.
\begin{remark}
	\Cref{eq:eupm_update} is of the same form as the update function used by \algmiff \cite{mifflin1993rapidly}.
\end{remark}

\begin{remark}
	We see the \acrshort{eupm} is qualitatively different from the \acrshort{supm} in that $\bracketstep{\eupmtag}$ depends only on the values of $\bracket$, and not on the function $f$ at all. 
	Therefore, like Golden Section, the \acrshort{eupm} is equally suited to smooth, non-smooth, unimodal and multi-modal functions and is scale invariant.
\end{remark}

\ifnum \condensed=0
\subsubsection{Theoretical Convergence of the \acrshort*{eupm}}
\else
\subsection{Theoretical Convergence of the \acrshort*{eupm}}
\fi
\label{sec:eupm_theoretical_bound}
In this section we derive a bound for the convergence rate of the \acrshort{eupm}. 
In particular, we find an integer $k$ and constant $C$ such that $b(\bracket_{i+k})/b(\bracket_i)\leq C$ holds for any objective function.
While there exist infinitely many possible ojective functions, there are only finitely possible values for $\bracket_{i+k}$ given $\bracket_i$.
This is because $\bracketstep{\eupmtag}$ does not depend on the objective function $f$ and thus the only effect that $f$ has on the \acrshort{eupm} is in the function $\bracketupdate{\eupmtag}$ (See \Cref{eq:UPM_iteration_update}) by determining which update function: $\bracketupdate{1}, \bracketupdate{2}, \bracketupdate{3}$ or $\bracketupdate{4}$  to apply, none of which depend on $f$. 

Therefore, there are only 4 possible values of $\bracket_{i+1}$ given $\bracket_i$: $\bracketupdate{j}\bracket$ for $j=1,2,3,4$. 
As a consequence of this, there are at most $4^k$ possible values of $\bracket_{i+k}$ given $\bracket_i$. 
Each possible value of $\bracket_{i+k}$ depends on the sequence of update functions which generated it. 

\begin{remark}
	In this section, we will refer to the update functions $\bracketupdate{1}, \bracketupdate{2}, \bracketupdate{3}$ and $\bracketupdate{4}$ extensively. 
	To simplify notation, we replace $\bracketupdate{j}$ with $\tilde{\bracketupdate{}}_j$, where $\tilde{\bracketupdate{}}_j(\bracket) := \bracketupdate{j}(\bracketstep{\eupmtag}(\bracket), \bracket)$.
	For the remainder of \Cref{sec:eupm}, $\bracketupdate{j}$ should be read as $\tilde{\bracketupdate{}}_j$.
\end{remark}

\begin{definition}
	\label{def:eupm_sequence}
	Let $\{ \bracket_i \}_{i=0}^n$ be a sequence of extended brackets generated by the \acrshort{eupm} such that $\bracket_{i+1} = \bracketupdate{\eupmtag}( \bracket_i)$ for each $i$. 
	Let $I= \{i_j\}_{j=0}^{n-1}\in \{1,2,3,4\}^n$ be the sequence of numbers such that $\bracket_{j+1} = \bracketupdate{i_{j}}(\bracket_j),\;\forall j \in \{0,1,\ldots,n-1\}$. 
	Then we call $I$ the sequence of \acrshort{eupm} iterations which generate $\{ \bracket_i \}_{i=0}^n$. 
	Moreover we write:
	\[ \bracketupdate{I} :=  \bracketupdate{i_n\circ}  \ldots  \bracketupdate{i_2\circ}   \bracketupdate{i_1},\; \text{where}\; I = \{i_j\}_{j=1}^n.  \]
\end{definition}
\begin{definition}
	\label{def:eupm_sequence_set}
	Define $\eupmsequencesetgen$ to be the set of all length $n$ sequences of \acrshort{eupm} iterations.
\end{definition}

\begin{definition}
	\label{def:function_projection}
	Define the $\funcprojection_k(f, \bracket_0)$ to be the sequence of the first $k$ iterations generated by the \acrshort{eupm} when applied to the function $f$ with the starting bracket $\bracket_0$.  
\end{definition}
\begin{lemma}
	\label{lem:func_projection_property}
	Let $f$ be a function, $\bracket \in \bracketset_f$ and $\eupmsequence \in \eupmsequencesetgen$ such that $\eupmsequence = \{\eupmsequence_1, \eupmsequence_2  \} = \funcprojection_n(f, \bracket)$, where $\eupmsequence_1 \in \eupmsequenceset{n-k}$ and $\eupmsequence_2 \in \eupmsequenceset{k}$. 
	Then it holds that: $\funcprojection_{k}(f, \bracketupdate{\eupmsequence_1}\bracket) = \eupmsequence_2$.
\end{lemma}
Using the notation from \Cref{def:eupm_sequence,def:function_projection}, it follows $\bracket_{n+k} = U_{\funcprojection_k(f,\bracket_n)}(\bracket_n)$. 
Since $\funcprojection_k(f,\bracket_n) \in \eupmsequenceset{k}$ for every function $f$, it follows that:
\[
\frac{b(\bracket_{n+k})}{b(\bracket_n)}\leq C\;\text{for any function f}\; \Leftrightarrow \frac{b(U_I \bracket_n)}{b(\bracket_n)}, \forall I \in \eupmsequenceset{k}.
\]

\begin{definition}
	\label{def:semi_bracket_set}
	Given $\eupmsequence  \in \eupmsequencesetgen$,
	we define: 
	\[\semibracketset_I:= \{\bracket \in \mathbb{R}^5|\;\exists f\;\text{such that}\;\bracket \in \bracketset_f,\; \funcprojection_n(f, \bracket)=I \}.\]
\end{definition}
\begin{lemma}
	\label{lem:semi_bracket_property}
	Let $\eupmsequence \in \eupmsequencesetgen$ be of the form $\eupmsequence = \{\eupmsequence_1, \eupmsequence_2 \}$, where $\eupmsequence_1 \in \eupmsequenceset{n-k}$ and $\eupmsequence_2 \in \eupmsequenceset{k}$.
	Then it holds that: $\semibracketset_\eupmsequence \subseteq \semibracketset_{\eupmsequence_1}$ and $\bracketupdate{\eupmsequence_1}(\semibracketset_\eupmsequence) \subseteq \semibracketset_{\eupmsequence_2}$.
\end{lemma}
\begin{proof}
If $\bracket \in \semibracketset_\eupmsequence$, then there exists a function $f$ such that $\eupmsequence=\funcprojection_n(f, \bracket)$. 
Since, $\eupmsequence = \{ \eupmsequence_1, \eupmsequence_2\}$, then $\eupmsequence_1 =\funcprojection_{n-k}(f, \bracket)$ by definition of $\funcprojection$.
Hence  $\bracket \in \semibracketset_{\eupmsequence_1}$.
Moreover, $\eupmsequence_2 =\funcprojection_{k}(f, \bracketupdate{\eupmsequence_1}\bracket)$ by \Cref{lem:func_projection_property}, which implies $\bracketupdate{\eupmsequence_1}\bracket \in \semibracketset_{\eupmsequence_2}$.
\end{proof}
We finally state our main result.
\begin{theorem}
	\label{thm:upm_limit_linear_convergence}
	The \acrshort{eupm} converges R-linearly. In particular:
	\[  \max \limits_{I \in \eupmsequenceset{5}} \max \limits_{\bracket \in \semibracketset_I} \frac{\bracketinnerdiam{ \bracketupdate{I} \bracket}}{\bracketinnerdiam{\bracket}} \leq \frac{1}{2}, \;\text{for any bracket}. \]
\end{theorem}

Given how long the proof of \Cref{thm:upm_limit_linear_convergence} is, we reserve an entire section for it. 
The proof may be found in \Cref{sec:proof_of_eupm_convergence}, although it may be skipped on first reading.



\Cref{thm:upm_limit_linear_convergence} provides a bound on the linear convergence rate for the \acrshort{eupm} albeit a weak one. 
In practice, we observe far faster convergence as will be demonstrated in \Cref{sec:univariate_numerics,sec:eupm_numerical_performance}.
However, the main justification for the \acrshort{eupm} is its robustness.
In \Cref{sec:upm}, we assumed that the \acrshort{supm} started with a bracket that was already sufficiently close to the solution. 
Once equipped with this, we might expect the \acrshort{supm} to converge quickly. 
The \acrshort{eupm} is useful for producing such a bracket in the first place. 
In \Cref{sec:dupm}, we present the final variation of the \acrshort{upm}, in which we combine the initial robustness of the \acrshort{eupm} with the desirable properties of the \acrshort{supm}.

\ifnum \condensed=0
	\subsubsection{Proof of \Cref{thm:upm_limit_linear_convergence}}

\else 
	\subsection{Proof of \Cref{thm:upm_limit_linear_convergence}}
\fi
\label{sec:proof_of_eupm_convergence}

We divide this proof into two parts, which are defined by \Cref{lem:minimal_set_works,lem:construct_minimal_set}.

\begin{lemma}
	\label{lem:minimal_set_works}
	It holds that $\bracketdiam{\bracketupdate{\eupmsequence}\bracket}/\bracketdiam{\bracket}\leq \tfrac{1}{2}\;\forall \eupmsequence \in \mathcal{A} $ where $\mathcal{A} = \left\{ 44\right.,$ $111,$ $143,$ $422,$ $414, $ $434,$ $1411,$ $1141,$ $1423, $ $4322, $ $4314, $ $\left.4114 \right\}$.   
\end{lemma}

\begin{lemma}
	\label{lem:construct_minimal_set}
	If $\bracketdiam{\bracketupdate{\eupmsequence}\bracket}/\bracketdiam{\bracket}\leq \tfrac{1}{2}\;\forall \eupmsequence \in \mathcal{A} $ then $  \bracketdiam{\bracketupdate{\eupmsequence}\bracket}/\bracketdiam{\bracket}\leq \tfrac{1}{2}\;\forall \eupmsequence \in \eupmsequencesetpart$, where $\mathcal{A}$ is as defined in \Cref{lem:minimal_set_works}.
\end{lemma}

Given that the proof of \Cref{lem:minimal_set_works} is entirely algebraic calculations, we leave it for \Cref{sec:eupm_calculations}.
In this section, we work towards proving \Cref{lem:construct_minimal_set}.
Ultimately, \Cref{thm:upm_limit_linear_convergence} follows from these two results.

\begin{definition}
	\label{substring}
	Let $\eupmsequence =\{i_j\}_{j=1}^n\in \eupmsequencesetgen$ be a length $n$ sequence of iterations of the \acrshort{eupm}. A sub-string of the sequence $\eupmsequence$ is a subsequence of the form $\eupmsubstring=\{i_j\}_{j = k_1}^{k_2}$ such that $ 1 \leq k_1 < k_2 \leq n$. Denote by $\eupmsubstringset(\eupmsequence)$ the set of sub-strings of the sequence $I$.
\end{definition}

The first step towards proving \Cref{lem:construct_minimal_set} is to show that if $\eupmsequence \in \eupmsequencesetpart$ is a length 5 sequence of iterations and $\eupmsubstring \in \eupmsubstringset(\eupmsequence)$ is a sub-string of $\eupmsequence$, then $\bracketdiam{\bracketupdate{\eupmsubstring}\bracket}/\bracketdiam{\bracket}\leq \tfrac{1}{2}$ implies $\bracketdiam{\bracketupdate{\eupmsequence}\bracket}/\bracketdiam{\bracket}\leq \tfrac{1}{2}$.


\begin{lemma}
	\label{lem:substring_property}
	Let $\bracket $ be an extended bracket of $f$. 
	Given $\eupmsequence\in \eupmsequencesetgen$,  and $J\in \eupmsubstringset(\eupmsequence)$, 
	it holds that:
	\[\max \limits_{\bracket\in \semibracketset_\eupmsequence} \frac{\bracketinnerdiam{\bracketupdate{\eupmsequence} \bracket}}{\bracketinnerdiam{\bracket}}\leq  \max \limits_{\bracket\in \semibracketset_\eupmsubstring } \frac{\bracketinnerdiam{\bracketupdate{\eupmsubstring} \bracket}}{\bracketinnerdiam{\bracket}} .\]
\end{lemma}
\ifnum \condensed=0
\begin{proof}
	Since $\eupmsubstring$ is a sub-string of $\eupmsequence$, it follows that we might write $\eupmsequence$ in the form: $\eupmsequence = \{\eupmsequence_1, \eupmsubstring, \eupmsequence_2\}$. 
	This implies that $\bracketupdate{\eupmsequence} = \bracketupdate{\eupmsequence_2} \bracketupdate{\eupmsubstring} \bracketupdate{\eupmsequence_1}$.
	Next we have:
	\begin{equation*} 
	\frac{\bracketinnerdiam{\bracketupdate{\eupmsequence} \bracket}}{\bracketinnerdiam{\bracket}}= \frac{\bracketinnerdiam{\bracketupdate{\eupmsequence_2} \bracketupdate{\eupmsubstring} \bracketupdate{\eupmsequence_1} \bracket}}{\bracketinnerdiam{\bracket}}= \frac{\bracketinnerdiam{\bracketupdate{\eupmsequence_2} \bracketupdate{\eupmsubstring} \bracketupdate{\eupmsequence_1} \bracket}}{\bracketinnerdiam{\bracketupdate{\eupmsubstring} \bracketupdate{\eupmsequence_1} \bracket}}\frac{\bracketinnerdiam{ \bracketupdate{\eupmsubstring} \bracketupdate{\eupmsequence_1} \bracket}}{\bracketinnerdiam{\bracketupdate{\eupmsequence_1} \bracket}}\frac{\bracketinnerdiam{ \bracketupdate{\eupmsequence_1} \bracket}}{\bracketinnerdiam{\bracket}}.
	\end{equation*}
	Now we use the fact that $\bracketinnerdiam{\bracketupdate{I_2} \bracketupdate{J} \bracketupdate{I_1} \bracket}/\bracketinnerdiam{\bracketupdate{J} \bracketupdate{I_1} \bracket} \leq 1$ and $\bracketinnerdiam{ \bracketupdate{I_1} \bracket}/\bracketinnerdiam{\bracket} \leq 1 $ to get:
	
	\begin{equation}
	\label{proof:eq:substring_property1}
	\max \limits_{\bracket\in \semibracketset_I } \frac{\bracketinnerdiam{\bracketupdate{\eupmsequence} \bracket}}{\bracketinnerdiam{\bracket}} \leq  \max \limits_{\bracket\in \semibracketset_I }\ \frac{\bracketinnerdiam{ \bracketupdate{J} \bracketupdate{I_1} \bracket}}{\bracketinnerdiam{\bracketupdate{I_1} \bracket}} \leq \max \limits_{\pmb s \in \semibracketset_{\{\eupmsubstring, \eupmsequence_2\}}}\ \frac{\bracketinnerdiam{ \bracketupdate{J} \pmb s}}{ \bracketinnerdiam{\pmb s}}\leq \max \limits_{\pmb s \in \semibracketset_{\eupmsubstring}}\ \frac{\bracketinnerdiam{ \bracketupdate{J} \pmb s}}{ \bracketinnerdiam{\pmb s}},
	\end{equation}	
	where the last two steps follow from \Cref{lem:semi_bracket_property}.
\end{proof}
\fi


In the next few lemmas, we show that not all sequences of $\{1,2,3,4\}$ are possible sequences of \acrshort{eupm} iterations given a particular starting bracket $\bracket$.


%

\begin{lemma}
	\label{lem:UPM_bigM_iter_type_sequence}
	Let $\eupmsequence= \{ i_j \}_{j=0}^{n-1}   \in \eupmsequencesetgen$ be a sequence of \acrshort{eupm} iterations.
	If $i_j = 1$ then $i_{j+1} \in \{1,4\}$. 
	Similarly, if $i_j = 2$ then $i_{j+1} \in \{2,3\}$. 
		
\end{lemma}
\begin{proof}
	If $i_j=1$, then by definition $\bracket_{j+1} = \bracketupdate{1}\bracket_{j}$.
	From this, we calculate:
	\begin{align*}
	\bracketstep{\eupmtag}(\bracket_{j+1}) - \extbracketmiddle_{j+1}&=\frac{\extbracketright{1}{\idxj+1}\extbracketright{2}{\idxj+1} - \extbracketleft{1}{\idxj+1} \extbracketleft{2}{\idxj+1}}{\extbracketright{1}{\idxj+1} + \extbracketright{2}{\idxj+1} - \extbracketleft{1}{\idxj+1} - \extbracketleft{2}{\idxj+1}} - \extbracketmiddle_{j+1},\\
	&= \frac{\extbracketright{1}{\idxj}\extbracketmiddle_{j} - \extbracketleft{1}{\idxj} \extbracketleft{2}{\idxj}}{\extbracketright{1}{\idxj} + \extbracketmiddle_{j} - \extbracketleft{1}{\idxj} - \extbracketleft{2}{\idxj}} - \frac{\extbracketright{1}{\idxj}\extbracketright{2}{\idxj} - \extbracketleft{1}{\idxj} \extbracketleft{2}{\idxj}}{\extbracketright{1}{\idxj} + \extbracketright{2}{\idxj} - \extbracketleft{1}{\idxj} - \extbracketleft{2}{\idxj}},\\
	&= \frac{-(\extbracketright{2}{\idxj} - \extbracketmiddle_j)(\extbracketright{1}{\idxj}-\extbracketleft{1}{\idxj})(\extbracketright{1}{\idxj} - \extbracketleft{2}{\idxj})}{(\extbracketmiddle_j + \extbracketright{1}{\idxj} - \extbracketleft{1}{\idxj} - \extbracketleft{2}{\idxj})(\extbracketright{1}{\idxj} + \extbracketright{2}{\idxj} - \extbracketleft{1}{\idxj} - \extbracketleft{2}{\idxj})}<0,
	\end{align*}
	since each $\bracket_j$ is an extended bracket.
	Given that $\bracketstep{\eupmtag}(\bracket_{j+1}) - \extbracketmiddle_{j+1}<0$, it follows from \Cref{eq:UPM_iteration_update} that $i_{j+1} \in \{1,4\}$.
	The equivalent argument applies for when $i_j = 2$.

\end{proof}


In the statement of \Cref{lem:UPM_bigM_iter_type_sequence}, one may notice a symmetry between the functions $\bracketupdate{1}$, $\bracketupdate{4}$ and $\bracketupdate{2}$, $\bracketupdate{3}$ respectively.
In the lemmas which follow, we derive a relation between $\bracketupdate{1},\bracketupdate{4}$ and $\bracketupdate{2},\bracketupdate{3}$ respectively such that if $\bracketinnerdiam{\bracketupdate{I}\bracket}/\bracketinnerdiam{\bracket}\leq \tfrac{1}{2}$ holds for some $I$, it will still hold after replacing $\bracketupdate{1}$ and $\bracketupdate{4}$ with $\bracketupdate{2}$ and $\bracketupdate{3}$ and vice versa.


\begin{definition}
	\label{def:Xmap}
	Let $\eupmsequence \in \eupmsequencesetgen$ be a sequence of \acrshort{eupm} iterations. 
	We define $\swappoints: \mathbb{R}^5 \rightarrow \mathbb{R}^5$ to be the function which satisfies $\swappoints(a,b,c,d,e) = -(e,d,c,b,a)$.
\end{definition}

We make use of the function $\swappoints$ because of the following useful properties.
\begin{lemma}
	\label{lem:Xmap_properties1}
	It holds that $\swappoints \swappoints \bracket = \bracket$, $\bracketstep{\eupmtag}(\swappoints \bracket) = -\bracketstep{\eupmtag}(\bracket)$, $ \bracketinnerdiam{\swappoints \bracket} = \bracketinnerdiam{\bracket}$, $\bracketupdate{1}\bracket = \swappoints \bracketupdate{2} \swappoints\bracket$, $\bracketupdate{2}\bracket = \swappoints \bracketupdate{1} \swappoints\bracket$, $\bracketupdate{3}\bracket = \swappoints \bracketupdate{4} \swappoints\bracket$ and $\bracketupdate{4}\bracket = \swappoints \bracketupdate{3} \swappoints \bracket$.
\end{lemma}
\ifnum \condensed=0
\begin{proof}
	\begin{align*}
	\swappoints \swappoints \bracket &= \swappoints \swappoints(\altbracketlefttwo,\altbracketleftone,\altbracketrightone,\altbracketrighttwo) = \swappoints (\altbracketrighttwo,\altbracketrightone,\altbracketleftone,p) = (\altbracketlefttwo,\altbracketleftone,\altbracketrightone,\altbracketrighttwo) = \bracket,\\
	\bracketdiam{\swappoints \bracket} &= \bracketdiam{\swappoints (\altbracketlefttwo,\altbracketleftone,\altbracketrightone,\altbracketrighttwo)} = \bracketdiam{\altbracketrighttwo,\altbracketrightone,\altbracketleftone,\altbracketlefttwo} = \altbracketleftone+\altbracketrightone = \bracketdiam{\bracket}.
	\end{align*}

	\begin{align*}
	\swappoints \bracketupdate{2} \swappoints \bracket &= \swappoints \bracketupdate{2}  \left( \begin{array}{c} \altbracketrighttwo\\\altbracketrightone\\\altbracketleftone\\\altbracketlefttwo\end{array}\right)=\swappoints \left( \begin{array}{c}\altbracketrightone\\  \frac{\altbracketleftone(\altbracketlefttwo+\altbracketleftone) - \altbracketrightone(\altbracketrightone+\altbracketrighttwo)}{\altbracketlefttwo+2(\altbracketleftone+\altbracketrightone)+\altbracketrighttwo}\\ \frac{(\altbracketleftone+\altbracketrightone)(\altbracketleftone+\altbracketrightone+\altbracketrighttwo)}{\altbracketlefttwo+2(\altbracketleftone+\altbracketrightone)+\altbracketrighttwo}  \\\altbracketlefttwo \end{array} \right) =  \left( \begin{array}{c}\altbracketlefttwo\\ \frac{(\altbracketleftone+\altbracketrightone)(\altbracketleftone+\altbracketrightone+\altbracketrighttwo)}{\altbracketlefttwo+2(\altbracketleftone+\altbracketrightone)+\altbracketrighttwo} \\  \frac{\altbracketleftone(\altbracketlefttwo+\altbracketleftone) - \altbracketrightone(\altbracketrightone+\altbracketrighttwo)}{\altbracketlefttwo+2(\altbracketleftone+\altbracketrightone)+\altbracketrighttwo} \\\altbracketrightone \end{array} \right)= \bracketupdate{1}\bracket,\\
	\swappoints \bracketupdate{3} \swappoints \bracket &= \swappoints \bracketupdate{3}  \left( \begin{array}{c} \altbracketrighttwo\\\altbracketrightone\\\altbracketleftone\\\altbracketlefttwo\end{array}\right)=\swappoints \left( \begin{array}{c}  \altbracketrighttwo\\\altbracketrightone\\  \frac{\altbracketleftone(\altbracketlefttwo+\altbracketleftone) - \altbracketrightone(\altbracketrightone+\altbracketrighttwo)}{\altbracketlefttwo+2(\altbracketleftone+\altbracketrightone)+\altbracketrighttwo} \\\frac{(\altbracketleftone+\altbracketrightone)(\altbracketleftone+\altbracketrightone+\altbracketrighttwo)}{\altbracketlefttwo+2(\altbracketleftone+\altbracketrightone)+\altbracketrighttwo}     \end{array} \right) =  \left( \begin{array}{c}  \frac{(\altbracketleftone+\altbracketrightone)(\altbracketleftone+\altbracketrightone+\altbracketrighttwo)}{\altbracketlefttwo+2(\altbracketleftone+\altbracketrightone)+\altbracketrighttwo}  \\  \frac{\altbracketleftone(\altbracketlefttwo+\altbracketleftone) - \altbracketrightone(\altbracketrightone+\altbracketrighttwo)}{\altbracketlefttwo+2(\altbracketleftone+\altbracketrightone)+\altbracketrighttwo}   \\  \altbracketrightone  \\  \altbracketrighttwo   \end{array} \right)= \bracketupdate{4}\bracket.
	\end{align*}
\end{proof}
\fi
Having shown that the function $\swappoints$ relates $\bracketupdate{1}$ to $\bracketupdate{4}$ and $\bracketupdate{4}$ to $\bracketupdate{3}$, we now show that this relation also holds for sequences of iterations.
\begin{definition}
	\label{def:Wmap}
	Let $I = \{i_j\}_{j=1}^n \in \eupmsequencesetgen$. Define $\swapupdates: \eupmsequencesetgen \rightarrow \eupmsequencesetgen$ to be the function:
	$\swapupdates(I) = \{ w(i_j) \}_{j=1}^n$, where $w(1)=2, w(2)=1, w(3)=4$ and $w(4)=3$.
\end{definition}

\begin{lemma}
	\label{lem:Xmap_func_projection}
	If $f$ and $g$ are functions such that $g(x) = f(-x)$ and $\bracket_0 \in \bracketset_f$, then $\swappoints \bracket \in \bracketset_g$.
	Moreover, if $\{\bracket^f_i\}_{i=0}^n$ and $\{\bracket^g_i\}_{i=0}^n$ are the sequences of extended brackets generated by the \acrshort{eupm} when applied to $f$ and $g$ starting from $\bracket_0$ and $\swappoints \bracket_0$ respectively, then it holds that $\bracket^f_j = \swappoints \bracket^g_j$ for every $j =0,1,\ldots, n$ and $\funcprojection_n(g, \swappoints\bracket) = \swapupdates(\funcprojection_n(f, \bracket))$.
\end{lemma}
\begin{proof}
	Note that $\swappoints \bracket \in \bracketset_g$ follows from \Cref{def:extended_bracket,def:Xmap}.
	Let $\{\bracket^f_i\}_{i=0}^n$ and $\{\bracket^g_i\}_{i=0}^n$ be the sequences of extend brackets generated by applying the \acrshort{eupm} with the starting brackets $\bracket_0$ and $\swappoints \bracket_0$ to $f$ and $g$ respectively. 
	Define:
	\begin{equation}
	\label{eq:bracketcheck}
	\bracketcheck(\bracket) = \bracketstep{\eupmtag}(\bracket) - \extbracketmiddle.
	\end{equation}
	
	Since $\bracketstep{\eupmtag}(\bracket) = - \bracketstep{\eupmtag}(\swappoints \bracket)$, it follows that $f(\bracketstep{\eupmtag}(\bracket_0)) = g(\bracketstep{\eupmtag}(\swappoints\bracket_0))$ and $\bracketcheck(\bracket_0) <0 \Leftrightarrow \bracketcheck(\swappoints \bracket_0) >0 $. 
	Given that $sign(\bracketcheck(\bracket))$ is precisely the condition in \Cref{eq:UPM_iteration_update} which determines whether to apply update functions $\bracketupdate{1}, \bracketupdate{4}$ or $\bracketupdate{2}, \bracketupdate{3}$, it follows that $\bracket^f_{1} = \bracketupdate{j} \bracket^f_0 =  \bracketupdate{j}\bracket_0 \Leftrightarrow \bracket^g_1 = \bracketupdate{W(j)}\bracket^g_0 = \bracketupdate{W(j)} \swappoints\bracket_0 $.
	
	Applying \Cref{lem:Xmap_properties1}, we note that $\bracketupdate{W(j)} \swappoints\bracket_0 = \swappoints \bracketupdate{j} \bracket_0$ which implies that $\bracket^g_1 = \swappoints \bracket^f_1$. 
	From this it follows inductively that $\bracket^f_j = \swappoints \bracket^g_j$ for all $j = 1,2, \ldots, n$.
	Moreover, if we write $\funcprojection_n(f, \bracket)$ and $\funcprojection_n(g, \swappoints\bracket)$ in the forms $\{i^f_j\}_{j=0}^{n-1}$ and $\{i^g_j\}_{j=0}^{n-1}$ respectively, it also follows inductively that $i^f_j = W(i^g_j)$ for each $j=0,1,\ldots, n-1$; in other words $\funcprojection_n(g, \swappoints\bracket) = \swapupdates(\funcprojection_n(f, \bracket))$.
\end{proof}

\begin{corollary}
	\label{cor:swappoints_property}
	If $\eupmsequence \in \eupmsequencesetgen$, then: $\swappoints(\semibracketset_\eupmsequence) = \semibracketset_{\swapupdates(\eupmsequence)}$ (recall \Cref{def:semi_bracket_set}).
\end{corollary}
\ifnum \condensed=0
\begin{proof}
	If $\bracket \in \semibracketset_\eupmsequence$, then there exists a function $f$ such that $\eupmsequence = \funcprojection_n(f, \bracket)$ and $\bracket \in \bracketset_f$.
	By defining $g(x):= f(-x)$, it follows from \Cref{def:extended_bracket} that: $\bracket \in \bracketset_f \Leftrightarrow \swappoints\bracket \in \bracketset_g$.
	By \Cref{lem:Xmap_func_projection}, $\funcprojection_n(g, \swappoints\bracket) = \swapupdates(\funcprojection_n(f, \bracket))$ 
\end{proof}
\fi

\begin{corollary}
	\label{cor:Xmap_properties2}
	It holds that $\bracketupdate{I}\bracket = \swappoints \bracketupdate{\swapupdates(I)} \swappoints \bracket,\; \forall I \in \eupmsequencesetgen$.
\end{corollary}
\begin{proof}
	Let $f$ be a function such that $\eupmsequence = \funcprojection(f, \bracket)$. 
	By \Cref{lem:Xmap_func_projection}, the function $g(x) = f(-x)$ has the property that  $\funcprojection_n(g, \swappoints\bracket) = \swapupdates(\funcprojection_n(f, \bracket))$ and $\bracketupdate{I}(\bracket)= \bracketupdate{\funcprojection_n(f, \bracket)}(\bracket) = \swappoints\bracketupdate{\funcprojection_n(g, \swappoints\bracket)}(\swappoints \bracket)= \swappoints\bracketupdate{W(I)}(\swappoints \bracket)$.
\end{proof}

The map $\swapupdates$ serves as a natural bijective between sequences of \acrshort{eupm} iterations which start with $1$ or $4$ and those which start with $2$ or $3$.
Dividing $\eupmsequencesetpart$ is useful due to the following result.
\begin{corollary}
	\label{cor:iteration_sequence_symmetry}
	Given $I\in \eupmsequencesetgen$ a sequence of \acrshort{eupm} iterations, it holds that:
	\[ \max \limits_{\bracket \in \semibracketset_\eupmsequence} \frac{\bracketdiam{\bracketupdate{\eupmsequence} \bracket}}{\bracketdiam{\bracket}} = \max \limits_{\bracket\in \semibracketset_{\swapupdates(\eupmsequence)}} \frac{\bracketdiam{\bracketupdate{\swapupdates(\eupmsequence)} \bracket}}{\bracketdiam{\bracket}} \]
\end{corollary}
\begin{proof}
	Using \Cref{lem:Xmap_properties1} and \Cref{cor:Xmap_properties2,cor:swappoints_property} it follows that:
	\[ \max \limits_{\bracket \in \semibracketset_\eupmsequence} \frac{\bracketdiam{\bracketupdate{\eupmsequence} \bracket}}{\bracketdiam{\bracket}} = \max \limits_{\bracket\in \semibracketset_\eupmsequence} \frac{\bracketdiam{\swappoints\bracketupdate{ \swapupdates(\eupmsequence)}\swappoints \bracket}}{\bracketdiam{\bracket}} = \max \limits_{\bracket\in \semibracketset_\eupmsequence} \frac{\bracketdiam{\bracketupdate{ \swapupdates(\eupmsequence)}\swappoints \bracket}}{\bracketdiam{\swappoints\bracket}}= \max \limits_{\pmb s\in \semibracketset_{\swapupdates(\eupmsequence)}} \frac{\bracketdiam{\bracketupdate{ \swapupdates(\eupmsequence)} \pmb s}}{\bracketdiam{\pmb s}}.\]

\end{proof}

We are now equipped to prove \Cref{lem:construct_minimal_set}.
For this proof, we write particular sequences of iterations such as $\{4,1,3\}$ simply as $413$. 

\begin{proof}[\Cref{lem:construct_minimal_set}]
	To prove this lemma, we list out elements of $\mathcal{I}_3, \mathcal{I}_4$ and $\eupmsequencesetpart$ which begin with either  $1$ or $4$ and mark them in the following way. 
	If a sequence is contained in $\mathcal{A}$, we overline it. 
	If a sequence contains a sub-string which is contained in $\mathcal{A}$, we underline the relevant sub-string.
	Finally, if a sequence contains a sub-string $J$ such that $\swapupdates(J) \in \mathcal{A}$, then we place square brackets around $J$. 
	We will show that all elements of $\eupmsequencesetpart$ can be marked in one of these three ways.
	First we list elements of $\mathcal{I}_3$ which start with either $1$ or $4$. 
	There is no need to list the sequences starting with $2$ or $3$ because each of those sequences is merely $\swapupdates(\eupmsubstring)$ for some $\eupmsubstring$ listed below.  
	\begin{align*}
	\overline{111}	&& 141	&& \overline{143} && 411  && \overline{422}  && 431  && 4[33] && \underline{44}1 && \underline{44}3 \\
	114	&& 142	&& 1\underline{44} && \overline{414}  && 423  && 432  && \overline{434} && \underline{44}2 && \underline{444}
	\end{align*}
	We see that only $7$ (or $14$ if you count the equivalent sequences starting with $2$ or $3$) of the sequences above are unmarked. 
	We take these seven sequences, and list all elements of $\mathcal{I}_4$ which begin with any of these $7$ sub-strings. 
	We neglect to list out elements of $\mathcal{I}_4$ which begin with a marked sub-string because these elements are themselves guaranteed to be marked .
	%
	\begin{align*}
	\overline{1141}	&& 1\underline{143}	&& \overline{1411}	&& 1\underline{422} && 4\underline{111}  && 4231  && 42[33] && 4[311] && \overline{4322}  \\
	1142	&&11\underline{44}	&& 1\underline{414}	&& \overline{1423} && \overline{4114}  && 4232  && 4[234] && \overline{4314} && 4[323] 
	\end{align*}
	Using the same procedure as earlier, we see only 3 sequences are left unmarked.
	We list out all elements of $\eupmsequencesetpart$ which begin with one of these three sequences and note that all of them are marked.
	\begin{align*}
	11\underline{422}	&& 42\underline{311}	&& 4[2322] &&	1\underline{1423}	&& 4[2314]	&& 42[323].
	\end{align*}
	If a sequence $J$ is marked in any of these ways, then it follows from \Cref{lem:substring_property} and \Cref{cor:iteration_sequence_symmetry} that $\bracketdiam{\bracketupdate{\eupmsequence}\bracket}/\bracketdiam{\bracket}\leq \tfrac{1}{2}\;\forall \eupmsequence \in \mathcal{A} $ implies  $\bracketdiam{\bracketupdate{J}\bracket}/\bracketdiam{\bracket}\leq \tfrac{1}{2}$.
\end{proof}
%

%

\ifnum \condensed=0
	\subsubsection{Numerical Performance}
	
	\label{sec:eupm_numerical_performance}
	In \Cref{sec:eupm_theoretical_bound} we bounded the convergence rate over 5 iterations of the EUPM.
	In practice we observe far superior convergence rates.
	From our observation, this seems to be because the extremal values of $\altbracket_0$ which yield the worst convergence rate over 5 iterations tend to yield best case performance over 6 iterations and so on.
	In this section, we show what performance might realistically be expected from the EUPM.
	
	In order to construct a suitable experiment, we first must ask what factors affect the convergence of the EUPM? 
	The answer to this question is the function to minimise $f$, along with the initial value of $\altbracket_0$. 
	However, the only effect that the function $f$ actually has is to determine whether $f(\tilde x) <f(x^M)$.
	This along with the current value of $\altbracket$ uniquely determines which update function will be used.
	Therefore, when testing the convergence rate of the EUPM, we do not test it on a range of functions, but rather for different binary sequences, where $1$'s mean $f(\tilde x) <f(x^M)$ and $0$'s the opposite.  
	
	The main advantage of this is that while the set $C([a,b])$ has infinite cardinality, the set of binary sequences of length $n$ is finite. 
	Therefore, given an initial bracket $\altbracket_0$, it is possible to determine how the EUPM will perform on literally any function.
	Given $\altbracket_0$ and $I \in \{0,1\}^n$, we plot the average convergence rate of the EUPM.
	
	\begin{figure}[h!]
		\centering

		\includegraphics[width=1\textwidth]{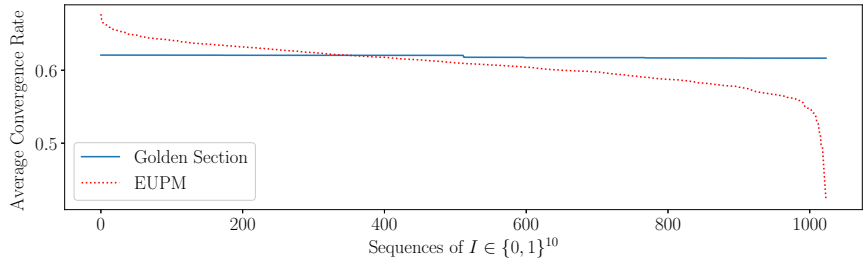}

		\caption{We compare the convergence rate of the EUPM with Golden Section over all possible sequences $I \in \{0,1\}^10$. These sequences are ordered such that the slowest ones are plotted first. For each sequences, we sample a large number of values for $\altbracket_0$ and average over these. }
		\label{fig:eupm_vs_gs}
	\end{figure}
	
	For each sequence $I \in \{0,1\}$, we sample a large number $\altbracket_0$ values, constrained such that $|\altbracket_0|_1 = 1$, and average over these.
	This leaves us with an average convergence rate for every possible sequence in $\{0,1\}$, and by extension for every possible function. 
	We choose $n=10$ and plot the resulting data corresponding to both the EUPM and Golden Section in \Cref{fig:eupm_vs_gs}.
	
	What we see is that Golden Section is the more conservative of the two. 
	It performs notably better in the worst case scenario even if notably worse in the best case scenario.
	Even though we see some convergence rates from the EUPM which are significantly worse than those observed in \Cref{tab:SUperformance,tab:NUperformance,tab:SMperformance}, they still are far better than the upper bound from \Cref{thm:upm_limit_linear_convergence}.
	\ifnum \condensed=0
	\subsection{Dynamic Underestimating Polynomial Method}
\else
	\section{Dynamic Underestimating Polynomial Method}
\fi
\label{sec:dupm}

In  \Cref{sec:upm}, we identified two main weaknesses of the SUPM: that of choosing a suitable $\upmconstant$, and that of ensuring that both $x^R_{1,i}$ and $x^L_{1,i}$ converge to $x^*$ fast enough.
In this section, we introduce a heuristic  which makes use of the features of the SUPM, while being equipped with a safeguarding option designed to avoid the situations where the SUPM fails. 
We call this heuristic the Dynamic UPM (DUPM), and denote the step and update functions for the DUPM by $\bracketstep{\dupmtag}$ and $\bracketupdate{\dupmtag}$ respectively.

The effect that $\upmconstant$ has on $\bracketstep{\supmtag}$ lies in how the model functions (see \Cref{eq:UPM_model_function_definitions}) are constructed.
We see that in this definition, $\upmconstant$ is multiplied to $\upmscaling(\bracket)$.
Since $\upmscaling(\bracket)\rightarrow 0$ should occur as the algorithm converges, a finite but excessively large value of $\upmconstant$ should not a problem ultimately. 
However, the DUPM may stall temporarily when $\upmconstant$ is too small.
Therefore, we construct the DUPM such that it will increase $\upmconstant$ when necessary, but never decrease it. 


In order for \Cref{thm:supm_convergence_ns_quality2} to apply, we need $\upmconstant$ to satisfy \Cref{lem:alpha_threshold}. 
While there is no way to check this, a necessary condition for \Cref{lem:alpha_threshold} to apply is that $f(x^M)\geq q^L(x^M)$ and $f(x^M)\geq q^R(x^M)$. 
Using a method similar to the proof of \Cref{lem:alpha_threshold}, we find that this is equivalent to:
\begin{equation}
\label{eq:alpha_threshold_2}
\upmconstant \geq \frac{1}{\upmscaling(\bracket)} \max \limits_{k \in \{L,R\}} \left(  f[x^k_1, x^k_2, x^k_3] - f[x^M, x^k_1, x^k_2]  \right).
\end{equation}
Therefore, at each iteration of the DUPM we set $\upmconstant_{i+1} \geq \max(\upmconstant_i, \upmconstant^*)$ where $\upmconstant^*$ is the smallest value $\upmconstant$ which satisfies \Cref{eq:alpha_threshold_2}.

Next we wish to force the DUPM to update both $x^L_{1,i}$ and $x^R_{1,i}$ regularly.
In order to do this, we must first equip the DUPM to recognize when this occurs. 
We append the three variables $u_1, u_2,u_3 $, whose purpose is record whether the update function $\bracketupdate{\dupmtag}$ updated $x^L$ or $x^R$ during each of the last three iterations, to the extended bracket $\bracket$. 
The update function $\bracketupdate{\dupmtag}$ is defined to be the natural extension of $\bracketupdate{\supmtag}$ (see \Cref{eq:UPM_iteration_update}) which includes $u_1, u_2$ and $u_3$:
\begin{equation}
\label{eq:UD}
\bracketupdate{\dupmtag}(\tilde x; \bracket, \upmconstant) = \left\{ \begin{array}{lcc}   
(\bracketupdate{1}(\tilde x;\bracket), R, u_1, u_2) & \text{if}&\tilde x<x^M\;\text{and}\;f(\tilde x) < f(x^M)\\
(\bracketupdate{2}(\tilde x;\bracket), L, u_1, u_2) & \text{if}& \tilde x>x^M\;\text{and}\;f(\tilde x) < f(x^M)\\
(\bracketupdate{3}(\tilde x;\bracket), R, u_1, u_2) & \text{if}& \tilde x>x^M\;\text{and}\;f(\tilde x) > f(x^M)\\	
(\bracketupdate{4}(\tilde x;\bracket), L, u_1, u_2) & \text{if}& \tilde x<x^M\;\text{and}\;f(\tilde x) > f(x^M)
\end{array}    \right.,
\end{equation}  
where the functions $\bracketupdate{1}, \ldots, \bracketupdate{4}$ are as defined in \Cref{eq:UPM_update_functions}.

If $u_1=u_2=u_3$, then this means that the DUPM has updated either $x^L_{1,i}$, or $x^R_{1,i}$ for each of the last three iterations.
When this occurs, we force the DUPM to take an EUPM step instead, $\bracketstep{\dupmtag}(\bracket;\upmconstant) = \bracketstep{\eupmtag}(\bracket)$.
This interference is similar to the fall back option used by Brent's method \cite{brent1976new,brentalgorithms}.
The choice to interfere after three iterations is based on practical experience, and not on any theoretical insight.

Finally, we define one more condition under which we interfere with the SUPM. 
From practical experience, we observe situations where either $q^L$ or $q^R$ is a convex function whose local minimum is returned by $\bracketstep{\supmtag}$, which result in slow convergence.
As a response to this, we ensure that $\upmconstant$ is sufficiently big that $\bracketstep{\dupmtag}$ returns a point of intersection between $q^L$ and $q^R$. 
More rigorously, we require $\upmconstant$ such that
\begin{equation}
\label{eq:ns_min_condition}
q^L\left(\bracketstep{\supmtag}(\bracket;\upmconstant); \bracket, \upmconstant \right)=q^R\left(\bracketstep{\supmtag}(\bracket;\upmconstant); \bracket, \upmconstant \right)
\end{equation}
holds, where $\bracketstep{\supmtag}$ is the step function for the SUPM (see \Cref{eq:UPM_gp}).

Define $\chi(\bracket)$ to be the the smallest value of $\upmconstant$ such that \Cref{eq:ns_min_condition} holds for all $\upmconstant \geq \chi(\bracket)$. 
Then at each iteration of the DUPM, we require $\alpha_{i+1} \geq \max(\alpha_i, \chi(\bracket))$.  

In order to compute $\chi(\bracket)$, note \Cref{eq:ns_min_condition} is satisfied whenever $q^L$ and $q^R$ are concave functions, provided that \Cref{eq:alpha_threshold_2} holds. 
Define $\upmconstant^+$ to be $\max_k f[x^k_1, x^k_2, x^k_3]/\upmscaling(\bracket)$, $k \in \{L,R\}$.
Then it follows that for all  $\upmconstant>\upmconstant^+$, \Cref{eq:ns_min_condition} is satisfied, implying that $\chi(\bracket) \leq \upmconstant^+$.

If at iteration $i$, $\upmconstant_i$ satisfies \Cref{eq:ns_min_condition}, then there is no need to compute $\chi(\bracket)$. 
Otherwise, we know that $\chi(\bracket)$ lies in the interval $(\upmconstant_i, \upmconstant^+]$.
Therefore, we can apply bisection to compute $\chi(\bracket)$ to any desired accuracy.

To define the DUPM rigorously, we must insert an additional line into \Cref{alg:bracketing_method_template}.
Before computing $\bracketstep{\dupmtag}$, we set
\[ \upmconstant_{i+1} = \max\left(\upmconstant_i, \chi(\bracket_i),   \frac{1}{\upmscaling(\bracket)} \max \limits_{k \in \{L,R\}} \left(  f[x^k_1, x^k_2, x^k_3] - f[x^M, x^k_1, x^k_2]  \right) \right). \]
Finally we define
\begin{equation}
\label{eq:gD}
\bracketstep{\dupmtag}(\bracket; \upmconstant):= \left\{ \begin{array}{cc}
\text{if}\;u_1=u_2=u_3 & \bracketstep{\eupmtag}(\bracket)\\
\text{else}& \bracketstep{\supmtag}(\bracket;  \upmconstant)
\end{array} \right..
\end{equation}

Ultimately the purpose of these interferences is to ensure that \Cref{thm:supm_convergence_ns_quality2} applies as soon as possible while forcing both $x^L_{1,i}$ and $x^R_{1,i}$ to be updated.

\ifnum \condensed=0
	\subsubsection{Numerical Performance}
	
	Having defined the DUPM, and justified parts of it based on observed results, all that remains is to show these results.
	We present a comparison of the DUPM vs Brent's method, Golden Section, the EUPM and SUPM in \Cref{tab:dupm_results}, using the same format as in \Cref{tab:supm_results}.

	\begin{table}
		\centering
		\pgfplotstabletypeset[	col sep = comma,
		/pgf/number format/precision=4,
		create on use/Functions/.style={
			create col/set list={$f^{SU}_1$,$f^{SU}_2$,$f^{SU}_3$,$f^{SU}_4$,$f^{SU}_5$,$f^{SU}_6$,$f^{SU}_7$,$f^{NU}_1$,$f^{NU}_2$,$f^{NU}_3$,$f^{NU}_4$,$f^{NU}_5$,$f^{SM}_1$,$f^{SM}_2$,$f^{SM}_3$,$f^{SM}_4$,$f^{SM}_5$,$f^{SM}_6$,$f^{SM}_7$}
		},
		columns/Functions/.style={string type},
		columns={Functions,0,0.1,1,10,EUPM,DUPM,Brent,Mifflin},
		display columns/1/.style={column name={$\alpha = 0$}},
		display columns/2/.style={column name={$\alpha = 0.1$}},
		display columns/3/.style={column name={$\alpha = 1$}},
		display columns/4/.style={column name={$\alpha = 10$}},
		display columns/5/.style={column name={EUPM}},
		display columns/6/.style={column name={DUPM}},
		display columns/7/.style={column name={Brent}},
		display columns/8/.style={column name={Mifflin}},
		every head row/.style={
			before row={\toprule}, 
			after row={\toprule}
		},
		every head row/.style={
			before row={
				\toprule
				&\multicolumn{4}{|c|}{SUPM with parameter value}&\multicolumn{4}{c|}{Other methods}
				\\
			},
			after row={ \toprule}
		},
		every last row/.style={after row=\bottomrule}, 
		]{CSVs/Univariatetest1000.csv}
		\caption{Average convergence rate of different algorithms on the set of smooth multimodal test functions.}
	\end{table}
	
	There are a few points to note from this table. 
	First and foremost, the DUPM performs at worst at the rate of Golden Section, but usually better. 
	Two of the three particular cases where the DUPM only matches Golden Section corresponds to the test-functions where $f''$ and $f^{(3)}$ also vanish at the minimiser. 
	On non-smooth functions, the DUPM outperforms Brent's method, except for $f^{NU}_5$, where Brent's method converges surprisingly fast.
	Finally, when comparing the DUPM with the SUPM, the DUPM converges at a rate between the SUPM for the two best values of $\upmconstant$, while it never fails to converge.
	
\fi
\else

\fi

\section{Numerical Results}
\label{sec:univariate_numerics}
Thus far, we have presented three algorithms: the SUPM, EUPM, and DUPM.
Of these, we have only bounded the theoretical convergence rate for EUPM. 
While \Cref{thm:supm_convergence_ns_quality2} is insightful, it does not yield any concrete bound. 
Moreover, while the DUPM has some nice properties which are designed to strengthen the SUPM, we have not offered any proof as to their effectiveness. 
Therefore, in order to justify these algorithms, we demonstrate their use on a variety of test functions, comparing them to Brent's method \cite{brentalgorithms,brent1976new} and \algmiff \cite{mifflin1990superlinear}.
Also, recall Golden Section which converges Q-linearly with rate $0.618$.

We consider three categories of test functions: smooth unimodal, non-smooth unimodal and smooth multimodal. 
To our knowledge there is not a commonly used set of univariate test functions. 
For our experiments, we have taken some test functions from Jamil et al \cite{jamil2013literature}, and  have added a few of our own, which are designed to be adversarial examples based on our observations. 
We scale each test function $f$ such that the minimum Lipschitz constant for $f'''$ lies between $0.8$ and $1$.
For the category of smooth unimodal, we use the following  test functions:
\begin{align}
	f^{SU}_1 &= -\frac{1}{\sqrt{e}}\exp\left( -\tfrac{1}{2} x^2 \right),& -1 \leq x \leq 1,\\
	f^{SU}_2 &= \frac{1}{24} x^4 ,& -1 \leq x \leq 1,\\
	f^{SU}_3 &= \frac{1}{11}\left(- \sin (2x - \tfrac{1}{2} \pi) - 3 \cos x - \tfrac{1}{2} x \right),& -2.5 \leq x \leq 3,\\
	f^{SU}_4 &= \frac{1}{2500}\left( \frac{x^2}{2} -\frac{ \cos \left(5 \pi x\right) }{25 \pi^2} - \frac{x \sin(5 \pi x)}{5 \pi } \right), &-10\leq x \leq 10,\\ 
	f^{SU}_5 &= -\frac{1}{250} \left(x^{\frac{2}{3}} + (1-x^2)^{\frac{1}{3}}\right),&0.1 \leq x \leq 0.9,\\
	f^{SU}_6 &= \frac{1}{6000} \left( e^x + \frac{1}{\sqrt{x}}\right) , & 0.1 \leq x \leq 3,\\
	f^{SU}_7 &= -\frac{1}{13}(16 x^2 - 24 x + 5) e^{-x},&1.3 \leq x \leq 3.9.
\end{align}

For each test function, we have an interval $[a,b]$ on which the function is intended to be used. 
We sample 4 points uniformly from both $[a, a+\tfrac{1}{5}(b-a)]$ and $[b - \tfrac{1}{2}(b-a), b]$.
After evaluating the function at these 8 points, we select the point where $f$ is minimised to be $x^M$; usually this will be one of the points closer to $\tfrac{1}{2}(b-a)$. 
Finally we order the remaining points, and select the 3 closest on the left and right to be $x^L_3, x^L_2, \ldots, x^R_3$. 
Selecting our initial points in this way yields an extended bracket.

Having generated the bracket, we run each algorithm and compute the average convergence rate over 1000 different random initialisations.
If an algorithm converged for a particular starting bracket, then its average convergence rate is by definition between 0 and 1, the smaller the better.
Iffor any random initialisation, an algorithm did not converge, we assign $\infty$ instead.
We present these results in \Cref{tab:SUperformance}.

\begin{table}
	\centering
	\pgfplotstabletypeset[	col sep = comma,
	/pgf/number format/precision=4,
	create on use/Functions/.style={
		create col/set list={$f^{SU}_1$,$f^{SU}_2$,$f^{SU}_3$,$f^{SU}_4$,$f^{SU}_5$,$f^{SU}_6$,$f^{SU}_7$}
	},
	columns/Functions/.style={string type},
	columns={Functions,0,0.1,1,10,EUPM,DUPM,Brent,Mifflin},
	display columns/1/.style={column name={$\alpha = 0$}},
	display columns/2/.style={column name={$\alpha = 0.1$}},
	display columns/3/.style={column name={$\alpha = 1$}},
	display columns/4/.style={column name={$\alpha = 10$}},
	display columns/5/.style={column name={EUPM}},
	display columns/6/.style={column name={DUPM}},
	display columns/7/.style={column name={Brent}},
	display columns/8/.style={column name={Mifflin}},
	every head row/.style={
		before row={\toprule}, 
		after row={\toprule}
	},
	every head row/.style={
		before row={
			\toprule
			&\multicolumn{4}{|c|}{SUPM with parameter value}&\multicolumn{4}{c|}{Other methods}
			\\
		},
		after row={ \toprule}
	},
	every last row/.style={after row=\bottomrule}, 
	]{SUtest1000.csv}
	\caption{Average convergence rate of different algorithms on the set of smooth unimodal test functions.}
	\label{tab:SUperformance}
\end{table}

First of all, we see that Brent's method consistently converges the fastest. 
This is not surprising, given that Brent's method is the only algorithm shown here which is specifically designed for smooth functions. 
While \algmiff is faster in select circumstances, it fails entirely on others.
This is also what we observe from the SUPM with $\upmconstant=0$. 

Looking at the SUPM's performance, we see that that the larger $\upmconstant$ is, the slower the SUPM will converge (assuming it converges in the first place).
The DUPM outperforms the SUPM with $\upmconstant=1$, but is not always faster than the SUPM with $\upmconstant=0.1$, which suggests that while \Cref{lem:alpha_threshold} is required for our theoretical results, satisfying it is by no means necessary.
Meanwhile, the EUPM is both consistent and slow, converging with an average rate comparable to Golden Section.

In short, the DUPM is not competitive with Brent's method for this class of test functions. 
However, it converge with an average rate which is better than for Golden Section, implying that it is not prohibitively slow.

Next we define our set of non-smooth unimodal test functions:
\begin{align}
	f^{NU}_1 &= -60000 \exp\left(-\frac{|x|}{50} \right),& -32 \leq x \leq 32,\\
	f^{NU}_2 &= \frac{1}{6} \max\left(\frac{1}{x+3}, \log(x) \right),& -2 \leq x \leq 10,\\
	f^{NU}_3 &=\frac{1}{24} \max \left(\frac{1}{x+3}, \frac{1}{(x-3)^2} \right),& -2 \leq x \leq 2,\\
	f^{NU}_4 &= \frac{1}{160}\max\left(\frac{1}{x+3}, \exp(x) \right), &-2\leq x \leq 5,\\ 
	f^{NU}_5 &= \frac{1}{150}\max\left( \exp(-x), \exp(x) \right),&-5 \leq x \leq 5.
\end{align}

The results for these functions are shown in \Cref{tab:NUperformance}.
These Results are particularly important as it was for this class of functions that the UPM was designed in the first place.
As we'd hope, we see the SUPM with $\upmconstant=1$ and the DUPM outperform Brent's method for every function.
Despite the fact that \algmiff is also designed for this context, it does not always successfully converge.
In particular, we observe that \algmiff struggles with non-convex functions.

Finally we define our smooth multimodal test functions.


%

\begin{table}
	\centering
	\pgfplotstabletypeset[	col sep = comma,
	/pgf/number format/precision=4,
	create on use/Functions/.style={
		create col/set list={$f^{NU}_1$,$f^{NU}_2$,$f^{NU}_3$,$f^{NU}_4$,$f^{NU}_5$}
	},
	columns/Functions/.style={string type},
	columns={Functions,0,0.1,1,10,EUPM,DUPM,Brent,Mifflin},
	display columns/1/.style={column name={$\alpha = 0$}},
	display columns/2/.style={column name={$\alpha = 0.1$}},
	display columns/3/.style={column name={$\alpha = 1$}},
	display columns/4/.style={column name={$\alpha = 10$}},
	display columns/5/.style={column name={EUPM}},
	display columns/6/.style={column name={DUPM}},
	display columns/7/.style={column name={Brent}},
	display columns/8/.style={column name={Mifflin}},
	every head row/.style={
		before row={\toprule}, 
		after row={\toprule}
	},
	every head row/.style={
		before row={
			\toprule
			&\multicolumn{4}{|c|}{SUPM with parameter value}&\multicolumn{4}{c|}{Other methods}
			\\
		},
		after row={ \toprule}
	},
	every last row/.style={after row=\bottomrule}, 
	]{NUtest1000.csv}
	\caption{Average convergence rate of different algorithms on the set of non-smooth unimodal test functions.}
	\label{tab:NUperformance}
\end{table}

%
\begin{align}
	f^{SM}_1 &= \frac{x^6}{300} \left(2 + \sin\left(\frac{1}{x} \right)\right),& -1 \leq x \leq 1,\\
	f^{SM}_2 &= -\frac{1}{80000}\sin\left(5 \pi x \right)^6 ,& -1 \leq x \leq 1,\\
	f^{SM}_3 &= -\frac{1}{250000} \sin\left(5 \pi (x^{\frac{3}{4}}-\frac{1}{20}) \right)^6,& 0.01 \leq x \leq 1,\\
	f^{SM}_4 &= \frac{1}{5}\left( \sin\left(\frac{16}{15} x - 1\right)+ \sin \left(\frac{16}{15} x - 1\right)^2 \right),  &-1\leq x \leq 1,\\ 
	f^{SM}_5 &= \frac{x^2}{4000} -\cos x + 1,&100 \leq x \leq 100,\\
	f^{SM}_6 &=\frac{1}{71} \left( \left(\log (x-2) \right)^2 +\left(\log (10-x) \right)^2 - x^{\frac{1}{5}}\right), & 2.5 \leq x \leq 9.5,\\
	f^{SM}_7&= \frac{1}{40}  \left(\sin(x) + \sin \left(\frac{10 x}{3}\right) + \log(x) + \frac{21}{25} x \right), & 0.5\leq x \leq 10.
\end{align}
It is important to note that theory for neither the UPM nor indeed Brent's method is suitable for multi-modal functions. 
However, both algorithms are equipped with features ensuring sufficient robustness in the general case that they converge (possibly slowly) to a locality where the function is unimodal.

For these test-functions, Brent's method is once again the best. 
While slower than Brent's method, the DUPM consistently performs better than Golden Section, which is an ideal result considering that these are not ideal circumstances for the DUPM.
Finally we see that \algmiff is entirely unsuitable for this type of problem, and as such consistently fails to converge in reasonable time. 
%
%

\begin{table}
	\centering
	\pgfplotstabletypeset[	col sep = comma,
	/pgf/number format/precision=4,
	create on use/Functions/.style={
		create col/set list={$f^{SM}_1$,$f^{SM}_2$,$f^{SM}_3$,$f^{SM}_4$,$f^{SM}_5$,$f^{SM}_6$,$f^{SM}_7$}
	},
	columns/Functions/.style={string type},
	columns={Functions,0,0.1,1,10,EUPM,DUPM,Brent,Mifflin},
	display columns/1/.style={column name={$\alpha = 0$}},
	display columns/2/.style={column name={$\alpha = 0.1$}},
	display columns/3/.style={column name={$\alpha = 1$}},
	display columns/4/.style={column name={$\alpha = 10$}},
	display columns/5/.style={column name={EUPM}},
	display columns/6/.style={column name={DUPM}},
	display columns/7/.style={column name={Brent}},
	display columns/8/.style={column name={Mifflin}},
	every head row/.style={
		before row={\toprule}, 
		after row={\toprule}
	},
	every head row/.style={
		before row={
			\toprule
			&\multicolumn{4}{|c|}{SUPM with parameter value}&\multicolumn{4}{c|}{Other methods}
			\\
		},
		after row={ \toprule}
	},
	every last row/.style={after row=\bottomrule}, 
	]{SMtest1000.csv}
	\caption{Average convergence rate of different algorithms on the set of smooth multimodal test functions.}
	\label{tab:SMperformance}
\end{table}

\section{Conclusion}
\label{sec:univariate_conclusion}
In this \ifnum \condensed=0 chapter \else paper\fi, we have constructed a univariate optimization algorithm for black box, piece-wise smooth functions. 
As seen in \Cref{sec:univariate_numerics}, this new method, the DUPM, converges both more robustly and often faster than existing methods for such problems.
Furthermore, while it is not the fastest algorithms for standard smooth test functions, it is not prohibitively slow either.

It must be acknowledged that the comparison between the UPM, \algmiff and Brent's method is not entirely fair given that they require a different number of points to start.
Therefore, in the context of a line-search, we expect to start with a method like Golden Section, and switch to the DUPM once we have accumulated enough points to form and extended bracket.
Regardless, the DUPM offers a univariate solver which performs well in most contexts and is therefore suitable for non-smooth functions.


\appendix

\section{ Proof of Lemma \ref{lem:minimal_set_works} }
\label{sec:eupm_calculations}

For the calculations in this section, we use the following change of variables:
$\altbracket = \left( \begin{array}{cccc} \altbracketlefttwo& \altbracketleftone&\altbracketrightone&\altbracketrighttwo\end{array}   \right)^{T} = \left( \begin{array}{cccc}x^L_1-x^L_2,& x^M- x^L_1,& x^R_1 - x^M,& x^R_2 - x^R_1\end{array}\right)^T.  $
Converting our calculations to being in terms of $\altbracket$ both reduces the number of variables, and simplifies the constraint $\bracket \in \semibracketset$ (recall \Cref{def:semi_bracket_set}) into $\altbracket \in \mathbb{R}^4_+$.

The algebra for the calculations which follow is tedious, and therefore we employed Mathematica to compute the values of $\bracketupdate{I}$, $I \in \mathcal{A}$ as well as changing the variables used.
\label{sec:eupm_calculations}
\begin{align*} 
\frac{\bracketdiam{ \bracketupdate{44} \bracket}}{\bracketdiam{\bracket}} &= \frac{(\altbracketleftone+\altbracketrightone)(\altbracketlefttwo+\altbracketleftone+\altbracketrightone)}{2(\altbracketleftone+\altbracketrightone)(\altbracketlefttwo+\altbracketleftone+\altbracketrightone) + (\altbracketleftone+\altbracketrightone)^2 +\altbracketrighttwo^2 + \altbracketrighttwo(\altbracketlefttwo+3\altbracketleftone+3\altbracketrightone) } < \frac{1}{2} ,\\
\frac{\bracketdiam{\bracketupdate{111} \bracket}}{\bracketdiam{\bracket}} &= \frac{\altbracketleftone}{\altbracketlefttwo+2\altbracketleftone+\altbracketrightone} < \frac{1}{2},\\
\frac{\bracketdiam{ \bracketupdate{143} \bracket}}{\bracketdiam{\bracket}} &= \frac{\altbracketleftone(\altbracketlefttwo+\altbracketleftone)(\altbracketlefttwo+\altbracketleftone+\altbracketrightone)}{ (\altbracketlefttwo+2\altbracketleftone+\altbracketrightone)(3\altbracketleftone^2 + 3\altbracketleftone\altbracketrightone + \altbracketrightone^2 + \altbracketlefttwo(2\altbracketleftone+\altbracketrightone))},\\
&< \frac{\altbracketleftone(\altbracketlefttwo+\altbracketleftone)}{3\altbracketleftone^2 + 3\altbracketleftone\altbracketrightone + \altbracketrightone^2 + \altbracketlefttwo(2\altbracketleftone+\altbracketrightone)} < \frac{1}{2},\\
\frac{\bracketdiam{ \bracketupdate{422} \bracket}}{\bracketdiam{\bracket}} &= \frac{(\altbracketleftone+\altbracketrightone)(\altbracketlefttwo+\altbracketleftone+\altbracketrightone)}{2(\altbracketleftone+\altbracketrightone)(\altbracketlefttwo+\altbracketleftone+\altbracketrightone) + (\altbracketleftone+\altbracketrightone)^2 +\altbracketrighttwo^2 + \altbracketrighttwo(\altbracketlefttwo+3\altbracketleftone+3\altbracketrightone) } < \frac{1}{2} ,\\
\frac{\bracketdiam{ \bracketupdate{1411} \bracket}}{\bracketdiam{\bracket}} &= \cfrac{1}{2+\left(\cfrac{ 4\altbracketleftone^3 + \altbracketlefttwo^2 \altbracketrightone+ 7\altbracketleftone^2 \altbracketrightone + 5\altbracketleftone\altbracketrightone^2 + \altbracketrightone^3+\altbracketlefttwo(\altbracketleftone+\altbracketrightone)(3\altbracketleftone+2\altbracketrightone) }{\altbracketleftone(\altbracketlefttwo+\altbracketleftone)(\altbracketlefttwo+\altbracketleftone+\altbracketrightone)}\right) }\leq \frac{1}{2},\\
\frac{\bracketdiam{ \bracketupdate{1141} \bracket}}{\bracketdiam{\bracket}} &= \frac{\altbracketleftone(\altbracketlefttwo+\altbracketleftone)(\altbracketlefttwo+\altbracketleftone+\altbracketrightone)}{(\altbracketlefttwo+2\altbracketleftone+\altbracketrightone)(\altbracketlefttwo^2 + 3\altbracketleftone^2 + \altbracketrightone(\altbracketrightone+\altbracketrighttwo) + 2\altbracketleftone(2\altbracketrightone+\altbracketrighttwo) + \altbracketlefttwo(3\altbracketleftone+2\altbracketrightone+\altbracketrighttwo))},  \\
&< \frac{\altbracketleftone(\altbracketlefttwo+\altbracketleftone)}{3\altbracketleftone(\altbracketlefttwo+\altbracketleftone)+ \altbracketlefttwo^2+\altbracketrightone(\altbracketrightone+\altbracketrighttwo)+2\altbracketleftone(2\altbracketrightone+\altbracketrighttwo)+2\altbracketrightone \altbracketlefttwo+\altbracketlefttwo\altbracketrighttwo } <\frac{1}{3},\\
\frac{\bracketdiam{ \bracketupdate{1423} \bracket}}{\bracketdiam{\bracket}} &
= \frac{(\altbracketlefttwo+\altbracketleftone+\altbracketrightone)(\altbracketlefttwo+2\altbracketleftone+\altbracketrightone)(\altbracketleftone(\altbracketlefttwo+\altbracketleftone )-\altbracketrightone(\altbracketrightone+\altbracketrighttwo))}{(\altbracketlefttwo+2\altbracketleftone+2\altbracketrightone+\altbracketrighttwo)} \;\times \\
&\frac{1}{(4\altbracketleftone(\altbracketleftone+\altbracketrightone)^2+\altbracketrightone^3 + \altbracketlefttwo^2(2\altbracketleftone+\altbracketrightone)+\altbracketleftone^2(\altbracketrighttwo-\altbracketrightone) + \altbracketlefttwo(6\altbracketleftone^2+2\altbracketrightone^2+\altbracketleftone(7\altbracketrightone+\altbracketrighttwo)))},\\&
<\frac{(\altbracketlefttwo+\altbracketleftone+\altbracketrightone)(\altbracketleftone(\altbracketlefttwo+\altbracketleftone )-\altbracketrightone(\altbracketrightone+\altbracketrighttwo))}{4\altbracketleftone(\altbracketleftone+\altbracketrightone)^2+\altbracketrightone^3 + \altbracketlefttwo^2(2\altbracketleftone+\altbracketrightone)+\altbracketleftone^2(\altbracketrighttwo-\altbracketrightone) + \altbracketlefttwo(6\altbracketleftone^2+2\altbracketrightone^2+\altbracketleftone(7\altbracketrightone+\altbracketrighttwo))},\\
&<\frac{(\altbracketlefttwo+\altbracketleftone+\altbracketrightone)(\altbracketleftone(\altbracketlefttwo+\altbracketleftone )-\altbracketrightone^2)}{4\altbracketleftone(\altbracketleftone+\altbracketrightone)^2+\altbracketrightone^3 + \altbracketlefttwo^2(2\altbracketleftone+\altbracketrightone)-\altbracketleftone^2 \altbracketrightone + \altbracketlefttwo(6\altbracketleftone^2+2\altbracketrightone^2+7\altbracketleftone\altbracketrightone)},\\
&<\frac{\altbracketleftone(\altbracketlefttwo+\altbracketleftone )-\altbracketrightone^2}{2\altbracketleftone(\altbracketlefttwo+\altbracketleftone)}<\frac{1}{2}.
\end{align*}

For the sequence $414$, we need to employ another piece of information. 
From \Cref{eq:UPM_iteration_update}, we see that $\bracketupdate{1}$ is only applied to $\bracket$ if $\bracketcheck(\bracket)<0$ (see \Cref{eq:bracketcheck}). 
This is equivalent to $\altbracketleftone(\altbracketlefttwo+\altbracketleftone)-\altbracketrightone(\altbracketrightone+\altbracketrighttwo)>0$ when written in terms of our alternative variables.
Proceeding with the calculation we find:
\begin{align*}
\frac{\bracketdiam{\bracketupdate{414} \bracket}}{\bracketdiam{\bracket}} &= \frac{1}{2+\left( \cfrac{\altbracketleftone^3 + \altbracketrightone^2(\altbracketlefttwo+\altbracketrightone) +\altbracketleftone^2(7\altbracketrightone+\altbracketrighttwo) +\altbracketleftone\altbracketrightone(3\altbracketlefttwo+7\altbracketrightone + 3\altbracketrighttwo) }{\altbracketleftone( \altbracketleftone(\altbracketlefttwo+\altbracketleftone)-\altbracketrightone(\altbracketrightone+\altbracketrighttwo))}\right)},
\end{align*}
which is bounded above by $\tfrac{1}{2}$ precisely when $\altbracketleftone(\altbracketlefttwo+\altbracketleftone)-\altbracketrightone(\altbracketrightone+\altbracketrighttwo)>0$.
Therefore $\bracketdiam{\bracketupdate{414} \bracket}/\bracketdiam{\bracket} <\tfrac{1}{2}$ as required.

There remain 4 elements of $\mathcal{A}$ for which we must establish a bound: $434$, $4322$, $4314$ and $4114$. 
As the previous change of variables does little to simplify the calculations for these subsequences, we instead use the following alternative change of variables: 
$(a,b,c,d)=\left( \begin{array}{cccc}x^R_1 - x^L_1,& x^R_1 - x^L_2,& x^R_2 - x^L_1,& \bracketcheck(\bracket)\end{array}\right).   $
For this set of variables, $\bracket \in \semibracketset$ is equivalent to $a,b,c>0$, $b>a$ and $c>a$. 
Using this transformation, we find:
\begin{align*}
\frac{\bracketdiam{ \bracketupdate{434} \bracket}}{\bracketdiam{\bracket}} &= \frac{b c^2 (b+c)(b+c-a)}{(a b+ b c + c^2)(a b^2 + 2 b^2 c + 3b c^2 + c^3)}<\frac{bc (b+c)}{a b^2 + 2 b^2 c + 3b c^2 + c^3} \\
&  < \frac{bc (b+c)}{2 b^2 c + 3b c^2 + c^3} < \frac{b (b+c)}{2 b^2  + 3b c + c^2}<\frac{b}{2b+c}<\frac{1}{2},\\
\frac{\bracketdiam{ \bracketupdate{4322} \bracket}}{\bracketdiam{\bracket}} &= \frac{b c^2 (b+c)(b+c-a)}{(a b+ b c + c^2)(a b^2 + 2 b^2 c + 3b c^2 + c^3)}<\frac{bc (b+c)}{a b^2 + 2 b^2 c + 3b c^2 + c^3} \\
&  < \frac{bc (b+c)}{2 b^2 c + 3b c^2 + c^3} < \frac{b (b+c)}{2 b^2  + 3b c + c^2}<\frac{b}{2b+c}<\frac{1}{2}.
\end{align*}
\begin{remark}
	From the calculations above, we see  that $\bracketdiam{ \bracketupdate{434} \bracket}/\bracketdiam{\bracket}=\bracketdiam{ \bracketupdate{4322} \bracket}/\bracketdiam{\bracket}$  and $\bracketdiam{ \bracketupdate{44} \bracket}/\bracketdiam{\bracket}=\bracketdiam{ \bracketupdate{422} \bracket}/\bracketdiam{\bracket}$.
	While this may imply another relation which we have not taken advantage of, we have not examined this further.
\end{remark}
Finally we have the sequences 4314 and 4114. 
Similar to what we did with the sequence $414$, we need to make use of additional information, specifically the fact that $\bracketupdate{4}$ is only applied if $\bracketcheck(\bracket)=d<0$.
In addition, $\bracketupdate{1}$ may only follow after $\bracketupdate{4}$ if $\bracketcheck(\bracketupdate{4}\bracket)<0$. 
This is equivalent to: $-d(ab+bc+c^2)>abc(b+c-a)$. 
From this it follows that:
\begin{align*}
\frac{\bracketdiam{\bracketupdate{4114} \bracket}}{\bracketdiam{\bracket}} &= \frac{a b c^2 (b+c)(b+c-a)}{(ab + bc+ c^2)(ac(b+c)^2-d(ab+bc+c^2)  ) },\\
&< \frac{a b c^2 (b+c)(b+c-a)}{(ab + bc+ c^2)( ac (b+c)^2 + abc(b+c-a)) },\\
&=\frac{b c(b+c)(b+c-a)}{(ab + bc+ c^2)( (2b+c)(b+c) -ab  )}<\frac{b c(b+c)^2}{(ab + bc+ c^2)( 2b+c)(b+c )},\\
&=\frac{b c(b+c)}{(ab + bc+ c^2)( 2b+c)}<\frac{b }{ 2b+c}<\frac{1}{2}.
\end{align*}

Similarly for the sequence $4314$, $\bracketupdate{3}$ may only follow after $\bracketupdate{4}$ if $\bracketcheck(\bracketupdate{4}\bracket)>0$. 
This is equivalent to: $-d(ab+bc+c^2)<abc(b+c-a)$.
First we compute 
\begin{align*}
\frac{\bracketdiam{\bracketupdate{4314} \bracket}}{\bracketdiam{\bracket}} &= \frac{-d(ac-d)(ab+bc+c^2)}{ a(b+c)(ac(b+c)^2-d(ab+bc+c^2)  ) }. 
\end{align*}
Next we observe that $\bracketcheck(\bracketupdate{4}\bracket)>0$ is equivalent to:
\begin{equation}
\label{eq:eupm_calc_1}
-d>\frac{abc(b+c-a)}{ab+bc+c^2} \Leftrightarrow ac - d >\frac{ac(b+c)^2}{ab+bc+c^2}.
\end{equation}
This implies that:
\begin{align*}
\frac{\bracketdiam{\bracketupdate{4314} \bracket}}{\bracketdiam{\bracket}} &< \frac{-dc (b+c)}{ ac(b+c)^2-d(ab+bc+c^2)   }. 
\end{align*}
Finally we that a sufficient condition for $\bracketdiam{\bracketupdate{4314} \bracket}/\bracketdiam{\bracket}<\tfrac{1}{2}$ is:
\begin{align*}
&-2dc (b+c) <  ac(b+c)^2-d(ab+bc+c^2)   ,\\
\Leftrightarrow & -d(-ab + bc+ c^2) < ac(b+c)^2.
\end{align*}
We see that this is true when we combine the requirement that $\bracketcheck(\bracketupdate{4}\bracket)>0$ must hold for the sequence $4314$ to occur with \Cref{eq:eupm_calc_1}.
In particular:
\[-d(-ab + bc+ c^2)<-d(ab + bc+ c^2)<(ac-d)(ab + bc+ c^2) <  ac(b+c)^2.\]
Therefore, $\bracketdiam{\bracketupdate{4314} \bracket}/\bracketdiam{\bracket}<\tfrac{1}{2}$ holds and moreover $\bracketdiam{\bracketupdate{\eupmsequence} \bracket}/\bracketdiam{\bracket}<\tfrac{1}{2}\;\forall \eupmsequence \in \eupmsequencesetpart$.

\section{Numerical Performance of the EUPM}


%
%

\bibliographystyle{spmpsci}
\bibliography{../../References/library}

%
%

\end{document}